\documentclass[a4paper]{elsarticle}
\usepackage[utf8]{inputenc}
\usepackage[T1]{fontenc}
\usepackage[english, francais]{babel}
\usepackage{amsmath}
\usepackage{amssymb}
\usepackage{amsthm}
\usepackage{color}
\usepackage{hyperref}

\newcommand{\md}{\mathrm{d}}

\newcommand{\R}{\mathbb{R}}
\newcommand{\C}{\mathbb{C}}
\newcommand{\N}{\mathbb{N}}

\renewcommand{\Re}{\mathrm{Re}}

\numberwithin{equation}{section}

\theoremstyle{definition}
\newtheorem{defi}{Definition}[section]

\theoremstyle{plain}
\newtheorem{prop}{Proposition}[section]
\newtheorem{theo}{Theorem}[section]
\newtheorem{coro}{Corollary}[section]  %Ces commandes s'utilisent aves des \begin{theo}
\newtheorem{lemme}{Lemma}[section]
\newtheorem{hypo}{Hypothesis}

\theoremstyle{remark}
\newtheorem{rmq}{Remark}[section]

\journal{J. Math. Pures Appl.}
\begin{document}

\begin{frontmatter}

\title{Quantitative Fattorini-Hautus test and minimal null control time for parabolic problems}

\author[LMB]{Farid \textsc{Ammar Khodja}}
\ead{fammark@univ-fcompte.fr}
\author[AMU]{Assia \textsc{Benabdallah}}
\ead{assia.benabdallah@univ-amu.fr}
\author[EDAN]{Manuel \textsc{Gonz\'alez-Burgos}\fnref{fn1}}
\ead{manoloburgos@us.es}
\author[AMU]{Morgan \textsc{Morancey}\corref{cor1}}
\ead{morgan.morancey@univ-amu.fr}

\fntext[fn1]{Supported by grant MTM2016-76990-P, Ministry of Economy and Competitiveness (Spain).}

\cortext[cor1]{Corresponding author}

\address[LMB]{Laboratoire de Math\'ematiques de Besan\c con, UMR 6623, Universit\'e de Franche-Comt\'e, 16 route de Gray, 25030 Besan\c con cedex, France.}
\address[AMU]{Aix Marseille Universit\'e, CNRS, Centrale Marseille, I2M, UMR 7373, 39 rue F. Joliot Curie 13453 Marseille, France.}
\address[EDAN]{Dpto.~Ecuaciones Diferenciales y An\'alisis Num\'erico and Instituto de Matem\'aticas de la Universidad de Sevilla (IMUS), Facultad de Matem\'aticas, Universidad de Sevilla, C/ Tarfia S/N, 41012 Sevilla, Spain.}

\begin{abstract}
\selectlanguage{english}
This paper investigates the link between the null controllability property for some abstract parabolic problems and an inequality that can be seen as a quantified Fattorini-Hautus test. 
Depending on the hypotheses made on the abstract setting considered we prove that this inequality either gives the exact minimal null control time or at least gives the qualitative property of existence of such a minimal time. 
We also prove that for many known examples of minimal time in the parabolic setting, this inequality recovers the value of this minimal time. 

\vskip 0.5\baselineskip

\selectlanguage{francais}
% Text of abstract in French
\noindent{\bf R\'esum\'e} \vskip 0.5\baselineskip \noindent
Dans cet article nous étudions le lien entre la contrôlabilité à zéro d'un problème parabolique abstrait et la validité d'une inégalité qui est une version quantifiée du test de Fattorini-Hautus.
Nous prouvons que cette inégalité permet de caractériser l'existence d'un temps minimal pour le problème de contrôlabilité à zéro et, selon les hypothèses considérées, d'obtenir la valeur de ce temps minimal.
Nous prouvons aussi que dans la plupart des exemples connus de problèmes paraboliques ayant un temps minimal de contrôle à zéro, cette inégalité est une conditions nécessaire et suffisante de contrôlabilité.
\end{abstract}
\selectlanguage{english}


\end{frontmatter}
\selectlanguage{english}

%\begin{abstract}

%\end{abstract}

%\maketitle

%%%%%%%%%%%%%%%%%%%%%%%%%%%%%%%%%%%%%%%%%%%%%%%%%%%%%%%%%%%%%%%%%%%%%%%%%%%%%%%%%%%%%%%%
\section{Introduction}
%%%%%%%%%%%%%%%%%%%%%%%%%%%%%%%%%%%%%%%%%%%%%%%%%%%%%%%%%%%%%%%%%%%%%%%%%%%%%%%%%%%%%%%%

%---------------------------------------------------------------------------------------
\subsection{Presentation of the minimal null control time problem}
%---------------------------------------------------------------------------------------

We consider an abstract control problem written in the form
\begin{equation} \label{SystControlAbstrait}
\left\{
\begin{aligned}
&y'(t) = A y(t) + B u(t), \quad &t& \geq 0,
\\
&y(0)=y_0.
\end{aligned}
\right. 
\end{equation}
In this equation, $u$ is a control that modifies the evolution of the state $y$.
In all what follows we consider that the operator $A$ is the infinitesimal generator of a $C_0$-semigroup on the Hilbert space $(H,\| \cdot \|)$ denoted by $t\mapsto e^{tA}$. 

The space of controls will be denoted by $(U,\|\cdot\|_{U})$. It is assumed to be a separable Hilbert space and is identified with its dual.
The inner products on $H$ and $U$ are respectively denoted by $\langle \cdot , \cdot \rangle$ and $\langle	 \cdot, \cdot \rangle_{U}$.
To take into account as many situations as possible we consider admissible control operators (see~\cite[Sec. 4.2]{TucsnakWeissBook} for instance), i.e. $B \in \mathcal{L}(U, D(A^*)')$ and its adjoint operator satisfies
\begin{equation} \label{DefAdmissibilite}
\exists K_T >0, \forall y \in D(A^*), \quad
\int_0^T \| B^* e^{tA^*} y\|_{U}^2 \md t \leq K_T \| y \|^2.
\end{equation}
This setting allows to consider, for any $y_0 \in H$ and any $u \in L^2(0,T;U)$, the unique solution $y \in C^0([0,T],H)$ of system~\eqref{SystControlAbstrait} defined by the following identity
\begin{equation} \label{DefSolution}
\langle y(t) , z^t \rangle - \langle y_0 , e^{t A^*} z^t \rangle = \int_0^t \langle u(\tau) , B^* e^{(t-\tau) A^*} z^t \rangle_U \md \tau, \quad \forall t \in [0, T],
\end{equation}
for any $z^t \in H$ (see for instance~\cite[Sec 2.3.1.]{CoronBook} or~\cite[Sec 4.2]{TucsnakWeissBook}).

The property under investigation in this article is null controllability. 
Recall that system~\eqref{SystControlAbstrait} is null controllable in time $T>0$ if for any $y_0 \in H$ there exists $u \in L^2(0,T;U)$ such that the associated solution of~\eqref{SystControlAbstrait} satisfies $y(T)=0$. 
We will also mention another notion of controllability: system~\eqref{SystControlAbstrait} is approximately controllable in time $T>0$ if for any $y_0, y_1 \in H$, for any $\varepsilon >0$, there exists $u \in L^2(0,T;U)$ such that the associated solution of~\eqref{SystControlAbstrait} satifies $\| y(T) - y_1 \| \leq \varepsilon$.

The following necessary condition for null controllability has been proved in~\cite[Prop 2.3]{DuyckaertsMiller12}. 
\begin{theo} \label{ThDuyckaertsMiller}
Assume that problem~\eqref{SystControlAbstrait} is null controllable in time $T$. Then, there exists $C_T >0$ such that for any $y \in D(A^*)$, for any $\lambda \in \C$ with $\Re(\lambda) >0$, one has
\begin{equation} \label{InegResolvante}
\| y \|^2 \leq C_T e^{2 T \Re(\lambda)} \left( \frac{\| (A^*+ \lambda) y\|^2}{\Re(\lambda)^2} + \frac{\|B^* y\|_{U}^2}{\Re(\lambda)} \right).
\end{equation}
\end{theo}
The aim of this article is to prove that this necessary condition is also sufficient for certain classes of parabolic control problems.

As detailed in~Sec.~\ref{Subsec:Biblio}, let us mention that this inequality can be seen as a quantified version of the Fattorini-Hautus test. 
In all this article we define
\begin{equation} \label{DefTempsMinimal}
\begin{aligned}
T^* = \inf &\bigg\{ T>0 \: ; \: \exists C_T >0 ; \forall y \in D(A^*), \forall \lambda \in \C \text{ with } \Re(\lambda) >0, 
\\
&\| y \|^2 \leq C_T e^{2 T \Re(\lambda)} \left( \frac{\| (A^*+ \lambda) y\|^2}{\Re(\lambda)^2} + \frac{\|B^* y\|_{U}^2}{\Re(\lambda)} \right) \bigg\}
\end{aligned}
\end{equation}
and $T^*=+\infty$ when the previous set is empty.

\medskip
The null controllability of parabolic partial differential equations has been widely studied since the pioneering work of \textsc{Fattorini} and \textsc{Russel}~\cite{FattoriniRussell71}. From the works of \textsc{Fursikov} and \textsc{Imanuvilov}~\cite{FursikovImanuvilov96} and \textsc{Lebeau} and \textsc{Robbiano}~\cite{LebeauRobbiano95}, it was commonly admitted that, in the context of parabolic partial differential equations, there is no restriction on the final time $T$ and no geometric restriction on the control domain (in case of internal or boundary control). 

But recently the study of particular examples highlighted the existence of a positive minimal time for null controllability or a geometric condition on the control domain. Actually, such an example was already provided in the $70$\textit{s} in~\cite{Dolecki73} but the full scope of this result was not understood at that time due to the particular pointwise control exerted. The more recent results concerning such a minimal time have been proved in contexts that have also been considered as specific, namely control of coupled parabolic equations~\cite{AKBDK05, AKBGBDT_condensation, AKBGBdT_JMAA16, Duprez_Tmin2017} or control of degenerate parabolic equations~\cite{BeauchardCannarsaGuglielmi, BHHR_Kolmogorov, BeauchardMillerMorancey_Tmin, BeauchardCannarsa17}. Though these three settings exhibit the same qualitative behavior, no precise link between them have been established so far.

Our aim in this article is to give an abstract framework encompassing those different frameworks to study the minimal null control time property. More precisely, we will relate this minimal time to the time $T^*$ defined in~\eqref{DefTempsMinimal}. 
We will emphasize that this minimal time can have different origins. It can be created by some localization of (generalised) eigenfunctions with respect to the observation operator $B^*$ (as in~\cite{Dolecki73, AKBGBdT_JMAA16, Duprez_Tmin2017, BeauchardCannarsaGuglielmi, BHHR_Kolmogorov, BeauchardMillerMorancey_Tmin, BeauchardCannarsa17}). This aspect is dealt with in Theorem~\ref{ThCondensationNulle}. But the minimal time can also be created by a condensation of the eigenvalues of the underlying operator as in~\cite{AKBDK05, AKBGBDT_condensation}. This aspect is dealt with in Theorem~\ref{ThCondensation}. 
In those two abstract settings the minimal null control time will exactly be given by $T^*$. 

We will also propose a more general setting (encompassing the two previous ones) to deal with the situation where the minimal time comes from both the localization of the eigenfunctions and the condensation of the spectrum. In this case (see Theorem~\ref{ThCondensationGeneral}), we will prove that the existence of such a minimal time is related to $T^*$ but the exact value of this minimal time will remain an open problem.

Finally there are still examples that will not fit into the different settings we study. For some of these examples (see Sec.~\ref{Sec:AutresCas}) we will still be able to prove that the minimal null control time is given by $T^*$. This analysis of particular examples will require to know a priori the value of the minimal time and thus, for the moment, the degenerate parabolic equations studied in~\cite{BeauchardCannarsaGuglielmi, BHHR_Kolmogorov, BeauchardCannarsa17} will remain out of the scope of this article.

%---------------------------------------------------------------------------------------
\subsection{Main results}
%---------------------------------------------------------------------------------------

Let us now define the setting that we will deal with in the next sections.
\begin{hypo} \label{HypoSpectrales}
Assume that the operator $-A^*$ admits a sequence of eigenvalues $\Lambda=(\lambda_k)_{k \in \N^*}$ such that 
\begin{equation} \label{ComportementSpectre}
\exists \delta >0, \: \Re(\lambda_k) \geq \delta |\lambda_k|, \: \forall k \in \N^*
\quad \text{and} \quad
\sum_{k=1}^{+\infty} \frac{1}{|\lambda_k|} <+\infty.
\end{equation}
For any $k \in \N^*$, we denote by $r_k = \mathrm{dim}( \mathrm{Ker}(A^*+\lambda_k))$ the geometric multiplicity of the eigenvalue $\lambda_k$ and assume that $\sup_{k \in \N^*} r_k < +\infty$. We denote by $(\varphi_{k,j})_{k\in \N^*, 1 \leq j \leq r_k}$ the associated sequence of normalised eigenfunctions and we assume that it forms a complete sequence in $H$, i.e.
\[
\Big( \langle \Phi , \varphi_{k,j} \rangle =0, \quad \forall k \in \N^*, \forall j \in \{1, \dots, r_k \} \Big) 
\: \Longrightarrow \: \Phi = 0.
\]
\end{hypo}

\begin{rmq}
As we will study parabolic problems, the hypothesis~\eqref{ComportementSpectre} mainly restricts our study to the one dimensional case. However, as presented in Sec.~\ref{Sec:AutresCas}, inequality~\eqref{InegResolvante} is also a sufficient condition for certain examples in higher dimensions such as the $2$D Grushin equation or a generalisation in any dimension of the academic example studied in~\cite{AKBDK05}.

Moreover, as underlined in Sec.~\ref{Subsec:Oscillateur}, hypothesis~\eqref{ComportementSpectre} is not a technical assumption due to the technics at stake in this article (namely the moment method) but is deeply related to the characterization we are studying.
\end{rmq}

\begin{rmq}
As it will be detailed in Remark~\ref{RmqDimensionValPropre}, one can weaken the assumption $\sup_{k \in \N^*} r_k < +\infty$ by authorizing a suitable growth.
\end{rmq}

\begin{rmq}
We will always assume that the sequence of eigenvalues is normally ordered that is
\begin{equation}\label{SuiteNormallyOrdered}
\left\{
\begin{aligned}
&\left\vert \lambda _{k}\right\vert \leq \left\vert \lambda_{k+1}\right\vert , \quad \forall k\geq 1, 
\\ 
&-\frac{\pi}{2} < \arg \left( \lambda_{k} \right) < \arg \left( \lambda_{k+1}\right) <\frac{\pi }{2} \hbox{ when } \left\vert \lambda _{k}\right\vert =\left\vert \lambda _{k+1}\right\vert.
\end{aligned}
\right.
\end{equation}
\end{rmq}

\begin{rmq}
Without loss of generality (through an orthonormalization process) we can always assume in Hypothesis~\ref{HypoSpectrales} that for any $k \in \N^*$,
\[
\langle \varphi_{k,i} , \varphi_{k,j} \rangle = \delta_{i,j}, \quad \forall i,j \in \{1, \dots, r_k\}.
\]
Here and in all this article $\delta_{i,j}$ denotes the Kronecker symbol, i.e. $\delta_{i,j}$ is equal to $1$ if $i=j$ and is equal to $0$ otherwise.
\end{rmq}

Before stating our results, let us recall the following definition of the condensation index of a complex sequence in this setting. 
\begin{defi} \label{DefiIndiceCondensation}
Let $\Lambda= (\lambda_k)_{k \in \N^*}$ be a complex normally ordered sequence satisfying~\eqref{ComportementSpectre}. Then the condensation index of this sequence is defined as
\[
\mathrm{c}(\Lambda) = \limsup_{k \to +\infty} \: \frac{- \ln |E'(\lambda_k)|}{\Re(\lambda_k)}
\]
where the function $E$ is defined by 
\begin{equation} \label{DefE} 
E: z \in \C \mapsto \prod_{k=1}^{+\infty} \left( 1 - \frac{z^2}{\lambda_k^2} \right). 
\end{equation}
\end{defi}
This definition can be found in~\cite{Shackell69}. The original (and of course equivalent) definition can be found in~\cite{Bernstein33}.
The first result of this article is the following theorem.
\begin{theo} \label{ThCondensationNulle}
Assume that $A^*$ satisfies Hypothesis~\ref{HypoSpectrales} and that the control operator $B$ is admissible. We also assume that the condensation index of the sequence $\Lambda=(\lambda_k)_{k\in \N^*}$ satisfies $\mathrm{c}(\Lambda)=0$. 
Let $T^*$ be defined by~\eqref{DefTempsMinimal}. Then:
\begin{itemize}
\item[$\bullet$] If $T^* >0$ and $T<T^*$, system~\eqref{SystControlAbstrait} is not null controllable in time $T$;
\item[$\bullet$] If $T^* < + \infty$ and $T>T^*$, system~\eqref{SystControlAbstrait} is null controllable in time $T$.
\end{itemize}
\end{theo}

\begin{rmq}
Notice that the negative result in time $T<T^*$ directly comes from Theorem~\ref{ThDuyckaertsMiller}.
\end{rmq}

\begin{rmq} \label{RmqGapCondensation}
It has to be highlighted that the condition $\mathrm{c}(\Lambda)=0$ is not a strong condition. It is more general than many situations previously considered when applying the moment method in the parabolic setting where a gap-like condition was assumed (see~\cite[Proposition 3.11]{AKBGBDT_Kalman}).

For example, the sequence $\Lambda = \left\{ k^2, k^2 +\frac{1}{k} ; k \in \N \right\}$ does not satisfy a gap-like condition whereas $\mathrm{c}(\Lambda)=0$ (see~\ref{AnnexeNoGap}). 
%This can be easily seen following the original definition using condensation groupings (see~\cite{Bernstein33, Shackell69}). Indeed this sequence admits only one condensation grouping given by $G_k = \{ k^2, k^2 + \frac{1}{k}\}$. As,
%\[
%\frac{-\ln \left( \frac{1}{k} \right)}{k^2}  \underset{k \to +\infty}{\longrightarrow} 0,
%\]
%the index of condensation of this grouping is, by definition, equal to $0$.
\end{rmq}

\begin{rmq} \label{RmqAnnonceJordan}
As it will be discussed in Sec.~\ref{Subsec:FurtherSyst_q} (on the particular case of system~\eqref{Syst_en_q}), the minimal time $T^*$ given by inequality~\eqref{DefTempsMinimal} can also turn out to be the minimal time for null controllability even if generalised eigenfunctions need to be taken into account to obtain a complete family. 
In the abstract setting we are considering, we were not able to prove (neither to disprove) this result. Instead, we prove in Proposition~\ref{PropJordan} a qualitative property: if a minimal time is needed for inequality~\eqref{InegResolvante} to hold, then a minimal time is also needed for null controllability. 
\end{rmq}

\medskip
Using inequality~\eqref{InegResolvante}, we may also characterize null controllability in certain settings where we allow the spectrum to condensate. First, we define the Bohr index (adapted from~\cite{Bernstein33} for a real-valued sequence) of a normally ordered sequence $\Lambda=(\lambda_k)_{k \in \N^*}$ by
\[
\mathrm{Bohr}(\Lambda) = \limsup\limits_{k \to +\infty} \: \frac{-\ln \inf\limits_{j \neq k }|\lambda_{k}-\lambda_j|}{\Re(\lambda_k)}.
\]
The Bohr index measures how two eigenvalues can get exponentially closer. The condensation index is in general a more subtle notion as it measures how packets of eigenvalues (with possibly increasing cardinality) get closer. Thus for any sequence $\Lambda$, it comes that $\mathrm{c}(\Lambda) \geq \mathrm{Bohr}(\Lambda)$. Under some additional assumptions, we prove that the equivalent of Theorem~\ref{ThCondensationNulle} still holds when $\mathrm{c}(\Lambda)=\mathrm{Bohr}(\Lambda)$, that is, that we allow eigenvalues to get closer two by two.
\begin{theo} \label{ThCondensation}
Assume that $A^*$ satisfies Hypothesis~\ref{HypoSpectrales} and that the eigenfunctions form a Riesz basis of $H$. Assume that the control operator $B$ is admissible and that for any $k \neq j \in \N^*$,
\begin{equation} \label{HypoStructureB}
\mathrm{Ker}(B^*) \cap \mathrm{Span}(\varphi_k, \varphi_{j}) \neq \{0\}
\end{equation}
and for any $\varepsilon >0$
\begin{equation} \label{ActionControleSuffisante}
\| B^* \varphi_k \|_{U} e^{\varepsilon \Re(\lambda_k)} \underset{k \to +\infty}{\longrightarrow} +\infty.
\end{equation}
Finally, assume that the condensation index of the sequence $\Lambda=(\lambda_k)_{k\in \N^*}$ satisfies $\mathrm{c}(\Lambda)=\mathrm{Bohr}(\Lambda)$. 

Let $T^*$ be defined by~\eqref{DefTempsMinimal}. Then:
\begin{itemize}
\item[$\bullet$] If $T^* >0$ and $T<T^*$, system~\eqref{SystControlAbstrait} is not null controllable in time $T$;
\item[$\bullet$] If $T^* <+\infty$ and $T>T^*$, system~\eqref{SystControlAbstrait} is null controllable in time $T$.
\end{itemize}
Moreover, in this setting it holds
\[
T^* = c(\Lambda).
\]
\end{theo}
In this result we obviously only consider the case $\mathrm{c}(\Lambda)>0$. Otherwise if $\mathrm{c}(\Lambda)=0$ the result stated in Theorem~\ref{ThCondensationNulle} is the same with less hypotheses.

\begin{rmq} 
In the case $T^* < +\infty$, assumption~\eqref{HypoStructureB} implies that the eigenvalues are geometrically simple. 
In fact, assume that for some $k \in \N^*$ we have $r_k \geq 2$. Then, it comes that there exist $\alpha_{k,1}, \alpha_{k,2} \in \C$ such that
\[
v_k := \alpha_{k,1} \varphi_{k,1} + \alpha_{k,2} \varphi_{k,2} \in \mathrm{Ker}(B^*) \backslash \{0 \}.
\]
But applying inequality~\eqref{InegResolvante} with $v=v_k$ and $\lambda= \lambda_k$ implies $v_k =0$. This is why we denoted simply the eigenfunctions by $\varphi_k$ instead of $\varphi_{k,1}$.
\end{rmq}

\begin{rmq}\label{RmqStructureB}
Notice that assumption~\eqref{HypoStructureB} is automatically satisfied when we consider a scalar control (for example a boundary control in dimension one), i.e. $U=\C$. However this assumption can also be satisfied by systems of coupled parabolic equations with a particular structure (see for instance~\eqref{SystIlya}). 
Actually, assumption~\eqref{HypoStructureB} does not need to be satisfied for any $k$ and any $j$. This is detailed in Remark~\ref{Rmq:HypoStructureB_precisee}.
\end{rmq}

\begin{rmq}
Some examples of problems satisfying these assumptions (as well as verifications of the hypothesis $\mathrm{Bohr}(\Lambda) = \mathrm{c}(\Lambda)$) are discussed in Sec.~\ref{Subsec:ExamplesCondensation}.
\end{rmq}

\begin{rmq}
It seems that hypothesis~\eqref{ActionControleSuffisante} is only technical; it allows us to study the effect of condensation of the spectrum when there is no localization of the eigenfunctions (with respect to $B^*$). However, if we do not assume this behavior of the observation operator $B ^*$ on the eigenfunctions, we would obtain only a qualitative result of existence of a minimal null controllability time belonging to $[T^*, T^*+c(\Lambda)]$ (see Remark~\ref{RmqTminBohr}). Such a qualitative result is stated below in Theorem~\ref{ThCondensationGeneral} in a more general setting.
\end{rmq}

%\begin{rmq} \label{RmqT0=condensationAnnonce}
%If every assumption of Theorem~\ref{ThCondensation} hold but we assume that $u$ is a scalar control (i.e. $U=\C$) instead of~\eqref{HypoStructureB}, then it comes that $T^* = c(\Lambda)$. This is proved in Remark~\ref{RmqT0=condensation}. 
%\end{rmq}

Finally, in the most general setting we study, we highlight a striking phenomenon: in the case where a positive minimal time is needed for inequality~\eqref{InegResolvante} to hold (i.e. $T^* \in (0,+\infty)$) then a positive minimal time is also needed for null controllability of system~\eqref{SystControlAbstrait}. More precisely, we prove the following result. 
\begin{theo} \label{ThCondensationGeneral}
Assume that $A^*$ satisfies Hypothesis~\ref{HypoSpectrales}. Assume that the control operator $B$ is admissible.
Let $T^*$ be defined by~\eqref{DefTempsMinimal}. Then, there exists $\widetilde{T} \in [T^* , T^*+\mathrm{c}(\Lambda)]$ such that
\begin{itemize}
\item[$\bullet$] if $\widetilde{T}>0$ and $T< \widetilde{T}$, system~\eqref{SystControlAbstrait} is not null controllable in time $T$;
\item[$\bullet$] if $\widetilde{T}<+\infty$ and $\operatorname{c}(\Lambda)<+\infty$, for $T> \widetilde{T}$, system~\eqref{SystControlAbstrait} is null controllable in time $T$.
\end{itemize}
\end{theo}

Notice that this setting is rather general: no restriction is imposed on the condensation of eigenvalues, we allow geometric multiplicity for the eigenvalues and the control operator is only assumed to be admissible. 

\begin{rmq}
In this general abstract setting we obtain the qualitative behaviour for the existence of a positive null control time. However the question of determining if its value is $T^*$ remains an open question, even in the explicit example given in Sec.~\ref{ExSystDeuxDiffusionsControlePonctuel}. Thus, in the settings of Theorems~\ref{ThCondensationNulle} and~\ref{ThCondensation}, we can conclude that the minimal null control time is exactly $T^*$ when there is either localization of eigenvectors or condensation of the spectrum but not when those phenomena appear at the same time. 
\end{rmq}

%---------------------------------------------------------------------------------------
\subsection{A brief review of previous results} \label{Subsec:Biblio}
%---------------------------------------------------------------------------------------

As mentioned, inequality~\eqref{InegResolvante} can be seen as a quantified Hautus test. The Hautus test refers to a well known controllability necessary and sufficient condition in finite dimension~\cite{Hautus69}. In the case where $H$ and $U$ have finite dimensions and the operators $A$ and $B$ are matrices then system~\eqref{SystControlAbstrait} is controllable in time $T>0$ if there are no eigenfunction of $A^*$ in the kernel of $B^*$. This can be rewritten as
\begin{equation} \label{TestHautus}
\mathrm{Ker}(A^*+\lambda) \cap \mathrm{Ker}(B^*) = \{ 0 \}, \qquad \forall \lambda \in \C.
\end{equation}
Actually, a more general result was proved by \textsc{Fattorini} in~\cite{Fattorini66} before the work of Hautus. Consider again that $H$ and $U$ are Hilbert spaces (actually the result can be extended to Banach spaces) and that $A$ generates a $C_0$-semigroup, has a compact resolvent and that the generalised eigenfunctions of $A^*$ form a complete sequence. Consider that the control operator $B$ is bounded. Then, the criterion~\eqref{TestHautus} is proved to be equivalent to some unique continuation property (in infinite time) that turns out to be equivalent to approximate controllability when the semigroup generated by $A$ is analytic. We also mention the work of \textsc{Olive}~\cite{Olive13} where the hypothesis on $B$ is weakened. 

Of course, inequality~\eqref{InegResolvante} implies~\eqref{TestHautus}. This is why we considered it as quantified Fattorini-Hautus test.

\medskip
The origin of this work is a conjecture of \textsc{Russell} and \textsc{Weiss}~\cite{RussellWeiss94}.
To describe it briefly, let $H$ and $U$ be two Hilbert spaces and $%
A:D\left( A\right) \subset H\rightarrow H$ the generator of a $C^{0}$-
semigroup (denoted by $e^{tA}$)  on $H.$ Consider the observability
operator $C:D\left( A^* \right) \rightarrow U$ which is assumed $A^*$-bounded and
admissible in infinite time: there are two positive constants $\sigma
_{\infty }$ and $\gamma _{\infty }$ such that for all $y \in D\left( A^*\right)
,$%
\[
\left\Vert Cy \right\Vert _{U}\leq \sigma _{\infty }\left\Vert A^*y\right\Vert
\]%
and%
\[
\int_{0}^{\infty }\left\Vert Ce^{tA^*}y\right\Vert _{U}~\md t\leq \gamma
_{\infty }\left\Vert y\right\Vert.
\]%
The cited authors proved that, with the previous assumptions, if the pair $%
\left( A^*,C\right) $ is exactly observable on $\left( 0,\infty \right) $,
namely%
\[
\exists m_{\infty }>0,~\int_{0}^{\infty }\left\Vert Ce^{tA^*}y\right\Vert
_{U}~\md t\geq m_{\infty }\left\Vert y\right\Vert, \quad \forall y\in D\left(
A^*\right) ,
\]%
then there exists $m>0$ such that for any complex number $\lambda $ with 
$\mathrm{Re}\left( \lambda \right) >0$ and any $y\in D\left( A^*\right) ,$%
\begin{equation}
\frac{\left\Vert \left( A^*+\lambda\right) y\right\Vert^{2}}{
\Re(\lambda)^{2}}+\frac{\left\Vert Cy\right\Vert
_{U}^{2}}{\Re(\lambda) }\geq m\left\Vert
y\right\Vert ^{2}.  \label{E}
\end{equation}%
The conjecture was then (see \cite[p. 2]{RussellWeiss94}): if (\ref{E}) holds
then the pair $\left( A^*,C\right) $ is exactly observable on $\left(
0,\infty \right) .$

In the same paper, they showed that if $A^*$ has a Riesz basis of
eigenfunctions and a strong extra condition on the eigenvalues is satisfied,
then the conjecture is true. But, some years later, \textsc{Jacob} and \textsc{Zwart}~\cite{JacobZwart04} proved by mean of a counterexample that the condition (\ref{E}%
) is not sufficient in general for exact observability on $\left( 0,\infty
\right) .$

Exact observability on $\left( 0,\infty \right) $ implies exact
observability in some time $T_{0}>0$ (see \cite[Proposition 6.5.2. p. 194]{TucsnakWeissBook}). For skew-adjoint operators $A,$ a lot of works with various assumptions show that the condition~\eqref{E} is sufficient for exact observability on $\left( 0,\infty \right) $. They also give an estimate of this finite time $T_{0}$ of observability (for more details, see \cite[Chapter 6 and the comments of this chapter]{TucsnakWeissBook} ).

%To our knowledge, only few papers deal with observability in finite time. 
The corresponding necessary condition for observability in time $T>0$ (see Theorem~\ref{ThDuyckaertsMiller}, relation~\eqref{InegResolvante} in this paper), due to \textsc{Duyckaerts} and \textsc{Miller} \cite{DuyckaertsMiller12}, asks again about the sufficiency of this condition for observability in time $T>0.$ In \cite{DuyckaertsMiller12}, the authors give, in a parabolic setting, an example showing that~\eqref{InegResolvante} can be true for any $T>0$ without the pair $\left( A^*,B^{\ast }\right) $ being final-time observable at \emph{any time} $T>0$. This example together with its implications in our context are given in Sec.~\ref{Subsec:Oscillateur}. In this example, the second part of~\eqref{ComportementSpectre} will not be satisfied. Let us also mention that these authors provide in the cited paper some positive results (sufficiency of condition~\eqref{InegResolvante} for any time $T>0$) for abstract parabolic equations. The remark that this inequality can also identify a positive minimal time of observability in a parabolic setting seems to be new.

\paragraph{Structure of the article}
\leavevmode\par

We end this introduction by giving the strategy of proof in the simpler case where the eigenvalues of $A^*$ are simple. Then, in Sec.~\ref{Sec:PreuveCondensationNulle} we prove the general case stated in Theorem~\ref{ThCondensationNulle}.
In this case, the minimal time will be inherited from a localization of the eigenfunctions with respect to the control operator. 
Section~\ref{Sec:PreuveCondensation} is dedicated to the proof of Theorem~\ref{ThCondensation} where we allow some condensation of the spectrum of the underlying operator. We also prove in this section the more general result of this article stated in Theorem~\ref{ThCondensationGeneral}.
Finally, in Sec.~\ref{Sec:AutresCas} we prove that the minimal null control time is given by $T^*$ in some examples that do not fit into the different abstract results and discuss the assumption~\eqref{ComportementSpectre}.

%---------------------------------------------------------------------------------------
\subsection{Strategy of proof}
\label{Subsec:Strategy}
%---------------------------------------------------------------------------------------

To detail the general idea used in this article, let us start by proving Theorem~\ref{ThCondensationNulle} in the simpler case of simple eigenvalues ($r_k=1$ for all $k \in \N^*$). For the sake of simplicity the associated eigenfunction $\varphi_{k,1}$ will simply be denoted by $\varphi_k$. The negative result follows directly from Theorem~\ref{ThDuyckaertsMiller}. Thus, assume that $T^* < +\infty$ and set $T>T^*$.

The proof of null controllability in time $T>T^*$ will rely on the moment method (initially developed in~\cite{FattoriniRussell71, FattoriniRussell74}) although the considered space of controls may not be of finite dimension.
From~\eqref{DefSolution}, it comes that for any $k \in \N^*$,
\begin{equation} \label{PbMomentAbstrait}
\langle y(T) , \varphi_k \rangle_{{ }} - \langle y_0 , e^{-\lambda_k T} \varphi_k \rangle_{{ }} 
= \int_0^T e^{-\lambda_k (T-t)} \langle u(t) , B^* \varphi_k \rangle_{U} \md t.
\end{equation}
Thus, as the family $(\varphi_k)_{k \in \N^*}$ is assumed to be complete, it comes that $y(T) = 0$ if and only if for every $k \in \N^*$, the control $u$ satisfies
\begin{equation} \label{PbMomentAbstrait2}
\int_0^T e^{-\lambda_k (T-t)} \langle u(t) , B^* \varphi_k \rangle_{U} \md t = - e^{-\lambda_k T} \langle y_0 , \varphi_k \rangle_{{ }}. 
\end{equation}

To apply the moment method let us recall the concept of biorthogonal family. 
\begin{defi}
Let $\sigma =(\sigma_k)_{k \in \N}$ be a complex sequence and $T>0$. We say that the family of functions $(q_k)_{k \in \N} \subset L^2(0,T;\C)$ is a biorthogonal family to the exponentials associated with $\sigma$ if for any $k, j \in \N$
\[
\int_0^T e^{-\sigma_j t} q_k(t) \md t = \delta_{k,j}.
\]
\end{defi}

The following result is proved in~\cite[Theorem 4.1]{AKBGBDT_condensation}.
\begin{prop} \label{FamilleBiorthogonale}
Let $T>0$ and let $\sigma =(\sigma_k)_{k \in \N}$ be a normally ordered complex sequence satisfying~\eqref{ComportementSpectre}. Then, there exists a biorthogonal family $(q_k)_{k \in \N}$ to the exponentials associated with $\sigma$ such that for any $\varepsilon > 0$ there exists $C_\varepsilon >0$ such that
\[
\| q_k \|_{L^2(0,T;\C)} \leq C_\varepsilon e^{ \Re(\sigma_k) (\mathrm{c}(\sigma) + \varepsilon)}, \quad \text{ for $k$ sufficiently large},
\]
where $\mathrm{c}(\sigma)$ is the condensation index of the sequence $\sigma$. 
\end{prop}

\begin{rmq}
In all what follows, $C_\varepsilon$ will always denote a constant that may vary from one line to another which only depends on the parameter $\varepsilon$.
\end{rmq}

We now generalise the strategy previously used for instance in~\cite{Lagnese83, CMV2017, AllonsiusBoyerMorancey16}, that is, we seek for a control $u$ solving~\eqref{PbMomentAbstrait2} in the following form
\begin{equation} \label{FormeControleAbstrait}
u(t) = -\sum_{k \in \N^*} \alpha_k  q_k(T-t) \frac{B^* \varphi_k}{\| B^* \varphi_k \|_{U}^2},
\end{equation}
where the complex coefficients $(\alpha_k)_{k \in \N^*}$  have to be determined. 

\begin{rmq}
Notice that when $u$ is a scalar control (i.e. $B^* : H \to \C$) then this is the classical form of controls used in the moment method. When $B$ is a localization operator and $u$ depends both on time and space, this is different from the usual strategy where one looks for a fixed space profile multiplied by a scalar control. In any other situation, this is a new way of applying the moment method.
\end{rmq}

First of all notice that in the previous form we divided by $\| B^* \varphi_k \|_{U}$. This is done only by convenience, as it could be incorporated in the constant $\alpha_k$, but we should justify that $\| B^* \varphi_k \|_{U} \neq 0$ for every $k \in \N^*$. 

Since we assumed $T^* < +\infty$ it comes that for any $\varepsilon>0$, for any $\lambda \in \C$ with $\Re(\lambda) >0$ and any $v \in D(A^*)$,
\[
\|v\|^2 \leq C_\varepsilon e^{2 \Re(\lambda) (T^* + \varepsilon)} \left( \frac{\|(A^*+\lambda)v\|^2}{\Re(\lambda)^2} + \frac{\|B^*v\|_{U}^2}{\Re(\lambda)} \right).
\]
The choice $\lambda= \lambda_k$ and $v= \varphi_k$ in the previous inequality implies that for all $k \in \N^*$,
\begin{equation} \label{EstimeeVectPropSimples}
\| B^* \varphi_k \|_{U} \geq C_\varepsilon \sqrt{\lambda_k} e^{-\mathrm{Re}(\lambda_k) (T^*+\varepsilon)} >0.
\end{equation}
Now, using the biorthogonality in the time variable, plugging $u$ given by~\eqref{FormeControleAbstrait} in the moment problem~\eqref{PbMomentAbstrait2} leads to 
\[
\alpha_k = - e^{-\lambda_k T} \langle y_0 , \varphi_k \rangle_{{ }},
\]
and thus a formal solution of~\eqref{PbMomentAbstrait} is given by 
\begin{equation} \label{SolutionPbMomentAbstrait}
u(t) = -\sum_{k \in \N^*} e^{-\lambda_k T} \frac{\langle y_0, \varphi_k \rangle_{{ }}}{\| B^* \varphi_k \|_{U}^2} q_k(T-t) B^* \varphi_k.
\end{equation}
The only remaining point is to prove that this series is indeed convergent in $L^2(0,T;U)$.
This comes directly from the estimate~\eqref{EstimeeVectPropSimples} and Proposition~\ref{FamilleBiorthogonale}. Indeed, for any $\varepsilon>0$, there exists $C_\varepsilon >0$ independent of $k$ such that
\begin{align*}
&\left\| \sum_{k \in \N^*} e^{-\lambda_k T} \frac{\langle y_0, \varphi_k \rangle_{{ }}}{\| B^* \varphi_k \|_{U}^2} q_k(T-\cdot) B^* \varphi_k \right\|_{L^2(0,T;U)} 
\\
&\leq \sum_{k \in \N^*} |\langle y_0, \varphi_k\rangle_{{ }}| e^{-\Re(\lambda_k) T} \frac{\| q_k \|_{L^2(0,T)}}{\|B^*\varphi_k \|_{U}} 
\\
&\leq C_\varepsilon \sum_{k \in \N^*} |\langle y_0, \varphi_k\rangle_{{ }}| e^{-\Re(\lambda_k) T} \frac{e^{\varepsilon \Re(\lambda_k)}}{\sqrt{\Re(\lambda_k)} e^{-\Re(\lambda_k)(T^* + \varepsilon)}}
\end{align*}
which converges for $\varepsilon$ sufficiently small, since $T>T^*$.

The strategy we develop in this article will be the same to prove Theorems~\ref{ThCondensationNulle} and~\ref{ThCondensation}. Namely we seek for a control in a suitable form and use inequality~\eqref{InegResolvante} to obtain sufficient estimates to prove that this control is well defined. 

\begin{rmq}
As a byproduct, we obtain that, in this setting, system~\eqref{SystControlAbstrait} is null controllable in a time $T>T^*$ if and only inequality~\eqref{EstimeeVectPropSimples} holds. 
\end{rmq}

%%%%%%%%%%%%%%%%%%%%%%%%%%%%%%%%%%%%%%%%%%%%%%%%%%%%%%%%%%%%%%%%%%%%%%%%%%%%%%%%%%%%%%%%
\section{Localization of eigenfunctions}
\label{Sec:PreuveCondensationNulle}
%%%%%%%%%%%%%%%%%%%%%%%%%%%%%%%%%%%%%%%%%%%%%%%%%%%%%%%%%%%%%%%%%%%%%%%%%%%%%%%%%%%%%%%%

In the first subsection we prove Theorem~\ref{ThCondensationNulle}. We then end this section by giving some examples in the literature of parabolic partial differential equations (with a zero condensation index for the underlying operator) exhibiting a minimal null control time given by~\eqref{DefTempsMinimal}.

%---------------------------------------------------------------------------------------
\subsection{Proof of Theorem~\ref{ThCondensationNulle}}
%---------------------------------------------------------------------------------------

The case $T<T^*$ follows  directly from Theorem~\ref{ThDuyckaertsMiller}. Thus we focus on the case $T^* <+\infty$ and we fix $T>T^*$.

As the moment problem is solved component by component, we distinguish two sets
\begin{equation*}
\Sigma_S = \{ k \in \N^* ; r_k =1 \}, 
\qquad
\Sigma_M = \{ k \in \N^* ; r_k \geq 2 \} = \N^* \backslash \Sigma_S,
\end{equation*}
where $r_k$ is the geometric multiplicity of the eigenvalue $\lambda_k$ defined in Hypothesis~\ref{HypoSpectrales}. For the sake of simplicity, when $r_k=1$, the associated eigenfunction $\varphi_{k,1}$ will simply be denoted by $\varphi_k$.

In Sec.~\ref{Subsec:Strategy}, we already proved Theorem~\ref{ThCondensationNulle} in the case where $\Sigma_S = \N^*$ thus we will assume in the following that $\Sigma_M \neq \varnothing$.

From a simple computation it comes that for any $k \in \N^*$, for any $i \in \{ 1, \dots, r_k \}$,
\begin{equation} \label{PbMomentMultiple}
\langle y(T) , \varphi_{k,i} \rangle_{{ }} - \langle y_0 , e^{-\lambda_k T} \varphi_{k,i} \rangle_{{ }} 
= \int_0^T e^{-\lambda_k (T-t)} \langle u(t) , B^* \varphi_{k,i} \rangle_U \md t.
\end{equation}
Thus as the family $(\varphi_{k,i})_{k \in \N^*, 1 \leq i \leq r_k}$ is assumed to be complete, it comes that $y(T) = 0$ if and only if for every $k \in \N^*$ and every $i \in \{ 1, \dots, r_k\}$, the control $u$ satisfies
\begin{equation} \label{PbMomentMultiple2}
\int_0^T e^{-\lambda_k (T-t)} \langle  u(t) , B^* \varphi_{k,i} \rangle_U \md t = - e^{-\lambda_k T} \langle y_0 , \varphi_{k,i} \rangle_{{ }}. 
\end{equation}
Using Proposition~\ref{FamilleBiorthogonale}, let $(q_k)_{k \in \N}$ be the biorthogonal family to the exponentials associated with $\Lambda=(\lambda_k)_{k \in \N^*}$. Following Sec.~\ref{Subsec:Strategy}, we define
\begin{equation} \label{FormeControleVPSimple}
u^1(t) = -\sum_{k \in \Sigma_S} e^{-\lambda_k T} \frac{\langle y_0, \varphi_k \rangle_{{ }}}{\| B^* \varphi_k \|_{U}^2} q_k(T-t) B^* \varphi_k.
\end{equation}
We already proved that $u^1$ is well defined in $L^2(0,T;U)$ and with this choice of control, the moment problem~\eqref{PbMomentMultiple2} is satisfied for every $k \in \Sigma_S$.

To deal with multiple eigenvalues we will need the following lemma whose proof is postponed at the end of this subsection.
\begin{lemme} \label{LemmeBiorthoEspace}
Let $(\mathcal{H},\langle \cdot, \cdot \rangle)$ be a Hilbert space and $r \in \N^*$. Let us consider $v = (v_1, \dots, v_r) \in \mathcal{H}^r$ a linearly independent family. We denote by
\[
G = \left( \langle v_i, v_j \rangle \right)_{1 \leq i,j \leq r}
\]
the associated Gram matrix. Then the matrix $G$ is invertible and the family $w=(w_1, \dots,w_r)$ defined by $w^t = G^{-1} v^t$ is biorthogonal to $v$, i.e.
\[
\langle v_i, w_j \rangle = \delta_{i,j}, \quad  \forall i,j \in \{1, \dots, r\}.
\]
Moreover, this biorthogonal family satisfies
\[
\| w \| \leq \frac{\sqrt{r}}{\sigma}, \quad \text{ where $\sigma^2$ is the smallest eigenvalue of $G$}.
\]
\end{lemme}
Let $k \in \Sigma_M$. To design our control we will apply this lemma to the family $(B^* \varphi_{k,1}, \dots, B^* \varphi_{k,r_k})$.

To this end we prove that it is linearly independent. Assume that there exist $c_1, \dots, c_{r_k} \in \C$ such that
\[
0 = c_1 B^* \varphi_{k,1} + \dots + c_{r_k} B^* \varphi_{k,r_k} = B^* \varphi,
\]
where $\varphi = c_1 \varphi_{k,1} + \dots + c_{r_k} \varphi_{k,r_k}$. Since $\varphi \in \mathrm{Ker}(A^* + \lambda_k)$, from~\eqref{DefTempsMinimal} and the assumption $T^* < +\infty$ it comes that $\varphi = 0$ which proves the linear independence.

\begin{rmq}
Notice that this construction would not be possible for a scalar control (i.e. $B^* : H \to \C$). This case is excluded in the presence of multiple eigenvalues thanks to the assumption $T^* < + \infty$. In this situation, system~\eqref{SystControlAbstrait} is not approximately controllable since one could combine two eigenfunctions (associated with the same eigenvalue) to design a solution which would not be observable. 
\end{rmq}

Using Lemma~\ref{LemmeBiorthoEspace}, we denote by $\Psi_k = (\Psi_{k,1}, \dots, \Psi_{k,r_k})$ such a biorthogonal family to $(B^* \varphi_{k,1}, \dots, B^* \varphi_{k,r_k})$.
We seek for a control $u^2$ in the form
\[
u^2(t) = \sum_{k \in \Sigma_M} q_k(T-t) \left(\sum_{1 \leq i \leq r_k} \alpha_{k,i} \Psi_{k,i} \right)
\]
such that with this choice of control, the moment problem~\eqref{PbMomentMultiple2} is satisfied for every $k \in \Sigma_M$. Plugging $u^2$ in~\eqref{PbMomentMultiple2} and using biorthogonality both in time and space variables leads to
\[
\alpha_{k,i} = -\langle y_0, \varphi_{k,i} \rangle_{{ }} e^{-\lambda_k T}.
\]
The only remaining point to prove is the convergence of the series defining $u^2$ in $L^2(0,T;U)$. 
The estimate on $q_k$ comes directly from Proposition~\ref{FamilleBiorthogonale}.

From Lemma~\ref{LemmeBiorthoEspace}, to estimate $\| \Psi_k \|_{U^{r_k}}$ one can simply give a lower bound on the smallest eigenvalue of the matrix
\[
G_k = \left( \langle B^* \varphi_{k,i} , B^* \varphi_{k,j} \rangle_U \right)_{1 \leq i, j \leq r_k}.
\]
For any $c = (c_1, \dots, c_{r_k})^t \in \C^{r_k}$ it comes that the vector defined by 
\[
\varphi = c_1 \varphi_{k,1} + \dots + c_{r_k} \varphi_{k,r_k}
\]
belongs to $\mathrm{Ker}(A^* + \lambda_k)$. From~\eqref{DefTempsMinimal} it comes that for any $\varepsilon >0$ there exists $C_\varepsilon>0$ such that for any $v \in \mathrm{Ker}(A^* + \lambda_k)$,
\[
\|  B^* v \|_{U}^2 \geq C_\varepsilon \Re(\lambda_k) e^{- 2\Re(\lambda_k) (T^* + \varepsilon)} \| v \|^2.
\]
On the other hand, a direct computation gives that $\langle G_k c, c \rangle_{\C^{r_k}} = \| B^* \varphi \|_{U}^2$. Gathering these two estimates it comes that for any $c \in \C^{r_k}$,
\begin{align*}
\langle G_k c, c \rangle_{\C^{r_k}} 
&\geq C_\varepsilon \Re(\lambda_k) e^{- 2\Re(\lambda_k) (T^* + \varepsilon)} \| c_1 \varphi_{k,1} + \dots + c_{r_k} \varphi_{k,r_k} \|^2
\\
&\geq C_\varepsilon \Re(\lambda_k) e^{- 2\Re(\lambda_k) (T^* + \varepsilon)} \left( c_1^2 + \dots + c_{r_k}^2 \right).
\end{align*}

Finally we obtain that for any $i \in \{1, \dots, r_k\}$,
\begin{equation} \label{EstimeeFamilleBiorthoEspace}
\| \Psi_{k,i} \|_U \leq C_\varepsilon \frac{\sqrt{r_k}}{\sqrt{\Re(\lambda_k)}} e^{\Re(\lambda_k) (T^* + \varepsilon)}
\end{equation}
Indeed, for any $\varepsilon>0$, there exists $C_\varepsilon >0$ such that
\begin{align*}
&\left\| \sum_{k \in \N^*} q_k(T-\cdot) \left( \sum_{1 \leq i \leq r_k} e^{-\lambda_k T} \langle y_0, \varphi_{k,i} \rangle_{{ }} \Psi_{k,i} \right) \right\|_{L^2(0,T;U)} 
\\
&\leq \| y_0 \| \sum_{k \in \N^*} e^{-\Re(\lambda_k) T} \| q_k \|_{L^2(0,T)}  \sum_{1 \leq i \leq r_k} \| \Psi_{k,i} \|_{U} 
\\
&\leq C_\varepsilon \|y_0\| \sup_{k \in \N^*} r_k^{\frac{3}{2}} \sum_{k \in \N^*} \frac{1}{\sqrt{\Re(\lambda_k)}} e^{-\Re(\lambda_k) T} e^{\varepsilon \Re(\lambda_k)}e^{\Re(\lambda_k)(T^* + \varepsilon)} 
\end{align*}
which converges since $T>T^*$ and $\sup\limits_{k \in \N^*} r_k < \infty$. 

Finally, notice that using biorthogonality with respect to the time variable $u^2$ does not affect the moment problem~\eqref{PbMomentMultiple2} for $k \in \Sigma_S$ and vice-versa.
Thus, setting $u=u^1+u^2$ ends the proof of Theorem~\ref{ThCondensationNulle}.

\hfill $\square$

\begin{rmq} \label{RmqDimensionValPropre}
Following closely the lines of the proof of Theorem~\ref{ThCondensationNulle} it is direct to remark that the hypothesis $\sup_{k \in \N^*} r_k < +\infty$ can be replaced by the following: for any $\varepsilon > 0$,
\[
r_k e^{-\varepsilon \Re(\lambda_k)} \underset{k \to + \infty}{\longrightarrow} 0.
\]
Nevertheless, we stick with $\sup_{k \in \N^*} r_k < +\infty$ for the sake of simplicity.
\end{rmq}

\begin{proof}[Proof of Lemma~\ref{LemmeBiorthoEspace}]
The matrix $G$ is hermitian. Since we assumed that the family $(v_1, \dots, v_r)$ is linearly independent, the Gram matrix $G$ is also positive definite. We define $w^t= G^{-1} v^t$. Let us notice that  $w$ is biorthogonal to $v$. Indeed, for any $i, j \in \{1, \dots, r\}$
\begin{align*}
\langle w_i , v_j \rangle
&= \langle \sum_{k=1}^r (G^{-1})_{i,k} v_k , v_j \rangle
\\
&= \sum_{k=1}^r (G^{-1})_{i,k} \langle v_k , v_j \rangle
\\
&=(G^{-1} G)_{i,j} = \delta_{i,j}.
\end{align*}
Finally the estimate on this biorthogonal family comes from the following identity
\[
r = \langle v, w \rangle = \langle Gw, w \rangle \geq \sigma^2  \| w \|^2.
\]
\end{proof}

%---------------------------------------------------------------------------------------
\subsection{Dealing with algebraically double eigenvalues}
%---------------------------------------------------------------------------------------

In this subsection we state and prove the qualitative result mentioned in Remark~\ref{RmqAnnonceJordan}.
We will use the following hypotheses.
\begin{hypo} \label{HypoSpectralesJordan}
Assume that the operator $-A^*$ admits a sequence of eigenvalues $\Lambda=(\lambda_k)_{k \in \N^*}$ satisfying~\eqref{ComportementSpectre}.

For any $k \in \N^*$ we assume that the algebraic multiplicity of the eigenvalue $\lambda_k$ is equal to $2$. We denote by $\varphi_{k,1}$ and $\varphi_{k,2}$, respectively, the associated normalised eigenfunction and generalised eigenfunction, i.e. 
\begin{equation} \label{DefFonctionPropreGeneralisee}
\left\{
\begin{aligned}
&-A^* \varphi_{k,1} = \lambda_k \varphi_{k,1},
\\
&-A^* \varphi_{k,2} = \lambda_k \varphi_{k,2} + \mu_k \varphi_{k,1},
\\
& \|\varphi_{k,1} \|= \|\varphi_{k,2}\| = 1.
\end{aligned}
\right.
\end{equation}
We assume that these eigenfunctions and generalised eigenfunctions form a complete sequence in $H$, i.e.
\[
\Big( \langle \Phi , \varphi_{k,j} \rangle =0, \quad \forall k \in \N^*, \forall j \in \{1, 2\} \Big) 
\: \Longrightarrow \: \Phi = 0.
\]
\end{hypo} 

\medskip
\begin{rmq}
Without loss of generality we will always assume that $\langle \varphi_{k,1}, \varphi_{k,2} \rangle = 0$. The orthonormalization process only affects the value of $\mu_k$.
\end{rmq}

\begin{prop} \label{PropJordan}
Assume that $A^*$ satisfies Hypothesis~\ref{HypoSpectralesJordan} and that $c(\Lambda)=0$.
Assume that $B$ is an admissible control operator and, for any $k \in \N^*$, one has
\begin{equation} \label{HypoStructureBJordan}
\mathrm{Ker}(B^*) \cap \mathrm{Span}(\varphi_{k,1}, \varphi_{k,2}) \neq \{0\}.
\end{equation}
Let $T^*$ be defined by~\eqref{DefTempsMinimal}. Then, there exists $\widetilde{T} \in [T^*, 2 T^*]$ such that
\begin{itemize}
\item[$\bullet$] If $T^* >0$ and $T< \widetilde{T}$, system~\eqref{SystControlAbstrait} is not null controllable in time $T$;
\item[$\bullet$] If $T^*<+\infty$ and $T> \widetilde{T}$, system~\eqref{SystControlAbstrait} is null controllable in time $T$.
\end{itemize}
\end{prop}

As for Theorem~\ref{ThCondensationGeneral}, we obtain a qualitative result of existence of minimal null control time but not its exact value.

\begin{rmq}
As for Theorem~\ref{ThCondensationNulle}, the proof relies on the moment method and the control is built separately for each eigenvalue. Thus, the same result also holds if the operator $-A^*$ admits a sequence of eigenvalues $\Lambda=(\lambda_k)_{k \in \N^*}$ satisfying~\eqref{ComportementSpectre} and the sequence $(\varphi_{k,j})_{k \in \N^*,1 \leq j \leq r_k}$ is complete in $H$ where, for every $k \in \N^*$,
\begin{itemize}
\item either we denote by $r_k = \mathrm{dim}( \mathrm{Ker}(A^*+\lambda_k))$ the geometric multiplicity of the eigenvalue $\lambda_k$ (assumed bounded) and by $(\varphi_{k,j})_{1 \leq j \leq r_k}$ the associated normalized eigenfunctions;
\item or the algebraic multiplicity of the eigenvalue $\lambda_k$ is equal to $2$ and we set $r_k = 2$ and denote by $\varphi_{k,1}$ and $\varphi_{k,2}$, respectively, the associated normalised eigenfunction and generalised eigenfunction as defined in~\eqref{DefFonctionPropreGeneralisee}.
\end{itemize}
Moreover as the main concern in the proof is the convergence of the series defining the control, if the latter has a finite occurrence, one recovers that $T^*$ is the minimal time for null controllability.
\end{rmq}

\begin{rmq}
The structural hypothesis made on $B^*$ is similar to~\eqref{HypoStructureB}. Thus, the comments of Remark~\ref{RmqStructureB} still hold in this setting.
\end{rmq}

\begin{proof}
As the construction of the control will be done separately for each eigenvalue, we can assume that $\mu_k \neq 0$ for any $k \in \N^*$: otherwise, for these frequencies, one simply follows the lines of Theorem~\ref{ThCondensationNulle}.

Now, as we are dealing with generalised eigenfunctions, the moment problem is different from~\eqref{PbMomentAbstrait2} or~\eqref{PbMomentMultiple2}. 
Let $y$ be the solution of~\eqref{SystControlAbstrait}. From~\eqref{DefSolution}, it comes that, for any $k \in \N^*$ and any $j \in \{1,2\}$, one has
\begin{equation}\label{PbMomentJordan}
\langle y(T), \varphi_{k,j} \rangle_{{ }} - \langle y_0, e^{T A^*} \varphi_{k,j} \rangle_{{ }}
= \int_0^T \langle Bu(t), e^{(T-t)A^*} \varphi_{k,j} \rangle_U.
\end{equation}
Using~\eqref{DefFonctionPropreGeneralisee}, we obtain 
\begin{equation}\label{SolutionVectPropresGeneralises}
\begin{aligned}
e^{(T-t)A^*} \varphi_{k,1} &= e^{-\lambda_k (T-t)} \varphi_{k,1}, 
\\
e^{(T-t)A^*} \varphi_{k,2} &= e^{-\lambda_k (T-t)} \varphi_{k,2} - (T-t) \mu_k e^{-\lambda_k(T-t)} \varphi_{k,1}.
\end{aligned}
\end{equation}
Using the assumption that the generalised eigenfunctions form a complete family, it comes that $y(T)=0$ if and only if, for any $k \in \N^*$ 
\begin{equation} \label{PbMomentJordan2}
\left\{
\begin{aligned}
& \int_0^T e^{-\lambda_k(T-t)} \langle u(t) , B^* \varphi_{k,1} \rangle_U \, \mathrm{d}t = 
- e^{-\lambda_k T} \langle y_0, \varphi_{k,1} \rangle_{{ }},
\\
& \int_0^T e^{-\lambda_k(T-t)} \langle u(t) , B^* \varphi_{k,2} \rangle_U \, \mathrm{d}t 
- \mu_k \int_0^T (T-t) e^{-\lambda_k(T-t)} \langle u(t) , B^* \varphi_{k,1} \rangle_U \, \mathrm{d}t 
\\
&= 
- e^{-\lambda_k T} \langle y_0 , \varphi_{k,2} - \mu_k T \varphi_{k,1} \rangle_{{ }}.
\end{aligned}
\right.
\end{equation}
The main difference with respect to the previous moment problem~\eqref{PbMomentAbstrait2} is the presence of the function $s \mapsto s e^{-\lambda_k s}$ in the second integral of the second equation. To deal with this, we will use an adapted biorthogonal family. 

We denote by $e_{k,1}$ and $e_{k,2}$ the $L^2(0,T;\C)$ functions defined by 
\[
e_{k,1}(t) = e^{-\lambda_k t}, 
\qquad
e_{k,2}(t) = t e^{-\lambda_k t}.
\]
To define and estimate a suitable biorthogonal family, we will use the following general proposition. Its proof is postponed to~\ref{AnnexeBiortho} and follows~\cite{AKBGBDT_Kalman, AKBGBDT_condensation}.
\begin{prop}\label{Prop:FamilleBiorthogonaleJordan}
Assume that $\Lambda=(\lambda_k)_{k \in \N^*}$ is a normally ordered sequence satisfying~\eqref{ComportementSpectre}.
Then, there exists a biorthogonal family $\{q_{k,1} , q_{k,2}\}_{k \in \N^*}$ to $\{e_{k,1}, e_{k,2}\}_{k\in \N^*}$ in $L^2(0,T;\C)$, i.e.
\[
\int_0^T e_{k,j}(t) q_{l,i}(t) \, \mathrm{d} t = \delta_{k,l} \delta_{j,i}, \quad \forall k,l \in \N^*, \forall j, i \in \{1,2\}.
\]
Moreover, we have the following estimate: for any $\varepsilon>0$, 
\begin{equation} \label{EstimeeBiorthoJordan}
\| q_{k,j} \|_{L^2(0,T)} \leq C_\varepsilon e^{(4c(\Lambda) + \varepsilon) \Re(\lambda_k)}, \quad 
\forall k \in \N^*, \forall j \in \{1,2\}.
\end{equation}
\end{prop}

Using this biorthogonal family we look for a control $u$ satisfying~\eqref{PbMomentJordan2} in the following form
\begin{equation} \label{FormeControleJordan}
u(t) = \sum_{k \in \N^*} \alpha_k q_{k,1}(T-t) \frac{B^* \varphi_{k,1}}{\| B^* \varphi_{k,1}\|_{U}^2} + 
\beta_k q_{k,2}(T-t) \frac{B^* \varphi_{k,1}}{\| B^* \varphi_{k,1}\|_{U}^2}.
\end{equation}
Thus, for $u$ to solve~\eqref{PbMomentJordan2}, the coefficients $\alpha_k$ and $\beta_k$ have to satisfy
\begin{equation} \label{CoeffMomentJordan}
\left\{
\begin{aligned}
& \alpha_k = - e^{-\lambda_k T} \langle y_0 , \varphi_{k,1} \rangle_{{ }},
\\
& \alpha_k \frac{\langle B^* \varphi_{k,2} , B^* \varphi_{k,1} \rangle_U}{\|B^* \varphi_{k,1}\|_{U}^2} - \mu_k \beta_k = - e^{-\lambda_k T} \langle y_0, \varphi_{k,2} - T \mu_k \varphi_{k,1} \rangle_{{ }}.
\end{aligned}
\right.
\end{equation}
From assumption~\eqref{HypoStructureBJordan} it comes that there exist $\alpha_{k,1}, \alpha_{k,2} \in \C$ such that
\[
\alpha_{k,1} \varphi_{k,1} + \alpha_{k,2} \varphi_{k,2} \in \mathrm{Ker}(B^*).
\]
As $T^* < +\infty$, inequality~\eqref{InegResolvante} implies that $B^* \varphi_{k,1} \neq 0$. Thus, there exists $\gamma_k \in \C$ such that $\gamma_k \varphi_{k,1} - \varphi_{k,2} \in \mathrm{Ker}(B^*)$ leading to
\[
B^* \varphi_{k,2} = \gamma_k B^* \varphi_{k,1}.
\]
Thus we obtain the following choice
\begin{equation} \label{CoeffMomentJordan2}
\left\{
\begin{aligned}
& \alpha_k = - e^{-\lambda_k T} \langle y_0 , \varphi_{k,1} \rangle_{{ }},
\\
& \beta_k = - e^{-\lambda_k T} \frac{\gamma_k}{\mu_k} \langle y_0 , \varphi_{k,1} \rangle_{{ }} - e^{-\lambda_k T} \frac{1}{\mu_k} \langle y_0, \varphi_{k,2} - T \mu_k \varphi_{k,1} \rangle_{{ }}.
\end{aligned}
\right.
\end{equation}
To conclude we need to prove that this defines a control $u \in L^2(0,T;U)$ when $T > 2T^*$. 

\medskip
Using inequality~\eqref{InegResolvante} it comes that for any $\varepsilon >0$,
\begin{equation} \label{EstimeeVectPropreJordan}
\| B^* \varphi_{k,1} \|_{U}^2 \geq C_\varepsilon \Re(\lambda_k) e^{-2\Re(\lambda_k)(T^*+\varepsilon)}.
\end{equation}
From the definition of $\gamma_k$, inequality~\eqref{InegResolvante} with $v = \gamma_k \varphi_{k,1} - \varphi_{k,2}$ implies for any $\varepsilon >0$
\begin{align*}
\| \mu_k \varphi_{k,1} \|^2 &\geq C_\varepsilon \Re(\lambda_k)^2 e^{-2 \Re(\lambda_k)(T^* + \varepsilon)} \| \gamma_k \varphi_{k,1} - \varphi_{k,2} \|^2 
\\
&= C_\varepsilon \Re(\lambda_k)^2 e^{-2 \Re(\lambda_k)(T^* + \varepsilon)} (1+\gamma_k^2).
\end{align*}
Thus we deduce that
\begin{equation} \label{Estimee_mu_kJordan}
|\mu_k| \geq C_\varepsilon \Re(\lambda_k) e^{- \Re(\lambda_k)(T^* + \varepsilon)}
\end{equation}
and
\begin{equation}\label{EstimeeGammaK}
\frac{|\gamma_k|}{|\mu_k|} \leq C_\varepsilon \frac{1}{\Re(\lambda_k)} e^{\Re(\lambda_k)(T^* + \varepsilon)} 
\end{equation}
\medskip
Using~\eqref{EstimeeBiorthoJordan} with the assumption $c(\Lambda)=0$ and~\eqref{EstimeeVectPropreJordan} it comes that for any $\varepsilon >0$,
\[
\left\| \alpha_k q_{k,1}(T-\cdot) \frac{B^* \varphi_{k,1}}{\|B^* \varphi_{k,1}\|_{U}^2} \right\|_{L^2(0,T;U)}
\leq C_\varepsilon \frac{|\langle y_0, \varphi_{k,1}\rangle_{{ }}|}{ \sqrt{\Re(\lambda_k)}} e^{-\Re(\lambda_k)T}  e^{\Re(\lambda_k)(T^* + 2\varepsilon)}.
\]
Similarly,
\begin{align*}
&\left\| \beta_k q_{k,2}(T-\cdot) \frac{B^* \varphi_{k,1}}{\|B^* \varphi_{k,1}\|_{U}^2} \right\|_{L^2(0,T;U)}
\\
&\leq C_\varepsilon e^{-\Re(\lambda_k)T} \frac{e^{\Re(\lambda_k)(T^* + 2\varepsilon)}}{\sqrt{\Re(\lambda_k)}}     
\left( \left(\frac{|\gamma_k|}{|\mu_k|}+T \right) |\langle y_0, \varphi_{k,1} \rangle_{{ }}| +
\frac{1}{|\mu_k|} |\langle y_0, \varphi_{k,2} \rangle_{{ }}| \right)
\\
&\leq C_\varepsilon \frac{1}{\Re(\lambda_k)^{3/2}}  e^{-\Re(\lambda_k)(T-2T^*-3\varepsilon)}  
\left( |\langle y_0, \varphi_{k,1} \rangle_{{ }}|  + |\langle y_0, \varphi_{k,2} \rangle_{{ }}|  \right).
\end{align*}
As $T> 2T^*$, the last two estimates end the proof of Proposition~\ref{PropJordan}.

\end{proof}

If we add another hypothesis on the behaviour of the control operator on the eigenfunctions then, we recover the exact value of the minimal time. More precisely, we get the following corollary.
\begin{coro} \label{CoroJordan}
Assume that the hypotheses of Proposition~\ref{PropJordan} hold. Assume moreover that for any $\varepsilon >0$,
\[
\| B^* \varphi_{k,1} \|_{U} e^{\varepsilon \Re(\lambda_k)}  \underset{k \to + \infty}{\longrightarrow} +\infty.
\]
Then, the conclusion of Proposition~\ref{PropJordan} hold with $\widetilde{T}= T^*$.
\end{coro}

\begin{proof}
The proof is exactly the same except that~\eqref{EstimeeVectPropreJordan} is now replaced by 
\[
\| B^* \varphi_{k,1} \|_{U} \geq C_\varepsilon e^{-\varepsilon \Re(\lambda_k)},
\]
for $k$ sufficiently large.
Following the lines leads to null controllability in time $T> T^*$.

\end{proof}

%---------------------------------------------------------------------------------------
\subsection{Examples of application} \label{Subsec:Examples}
%---------------------------------------------------------------------------------------

\subsubsection{Pointwise controllability of a one dimensional heat equation}
\label{ExempleDolecki}

Let us consider the system
\begin{equation} \label{SystDolecki}
\left\{
\begin{aligned}
& \partial_t y = \partial_{xx} y + \delta_{x_0}u(t), \quad &(t,x) \in (0,T) \times (0,1),
\\
& y(t,0) = y(t,1) = 0,    &t \in (0,T),
\\
& y(0,x) = y_0(x),        &x \in (0,1).
\end{aligned}
\right.
\end{equation}
The eigenfunctions and eigenvalues of the underlying operator are denoted by
\begin{equation} \label{EltsPropresClassiques}
\varphi_k(x) = \sqrt{2} \sin(k \pi x), \qquad
\lambda_k = k^2 \pi^2.
\end{equation}
In this case, the condensation index of the sequence $\Lambda=\{\lambda_k ; k \in \N\}$ is equal to $0$ (see for instance Remark~\ref{RmqGapCondensation}).
It is proved in~\cite{Dolecki73} that there is a minimal time for null controllability of~\eqref{SystDolecki} in $H=L^2(0,1)$ with controls in $L^2(0,T;\R)$ given by
\[
T_{min} = \limsup\limits_{k \to +\infty} \: \frac{- \ln |\varphi_k(x_0)|}{\lambda_k}.
\]
This is historically the first example of a minimal null control time in the parabolic framework.

The abstract setting to fit into the formalism of system~\eqref{SystControlAbstrait} is given for example in~\cite[Sec. 6.1]{AKBGBDT_condensation}. With this abstract setting it is straightforward to prove that this example satisfies the assumption of Theorem~\ref{ThCondensationNulle}. Thus we deduce that $T^* = T_{min}$.

\subsubsection{A parabolic cascade system with a variable coupling coefficient: a particular case of internal control}
\label{ExempleSyst_q_CasParticulier}

Let $q \in L^\infty(0,1;\R)$. Let us consider the following system studied in~\cite{AKBGBdT_JMAA16}
\begin{equation} \label{Syst_en_q}
\left\{
\begin{aligned}
& \partial_t y + \begin{pmatrix}  -\partial_{xx} & 0 \\ 0 & -\partial_{xx} \end{pmatrix} y + \begin{pmatrix} 0 & q(x) \\ 0 & 0 \end{pmatrix} y  = \begin{pmatrix} 0 \\ u \mathbf{1}_\omega \end{pmatrix},  &(t,x) \in (0,T)\times(0,1),
\\
& y(t,0) = y(t,1) = \begin{pmatrix} 0 \\ 0 \end{pmatrix},   &t \in (0,T),
\\
&y(0,x) = y_0(x),    &x \in (0,1),
\end{aligned}
\right.
\end{equation}
with $\omega = (a,b) \subset (0,1)$ and $\mathrm{Supp}(q) \cap \omega = \emptyset$.
For any $k  \in \N^*$ let
\begin{equation} \label{Def_Ik(q)}
I_k(q) = \int_0^1 q(x) \varphi_k(x)^2 \md x,
\qquad
I_{1,k}(q) = \int_0^a q(x) \varphi_k(x)^2 \md x. 
\end{equation}

To fit into our abstract setting, we set $H=L^2(0,1; \R^2)$, $U=L^2(\omega)$ and 
\begin{gather*}
A = \frac{\md^2}{\md x^2} \mathrm{Id} -  \begin{pmatrix} 0 & q(x) \\ 0 & 0 \end{pmatrix},
\\
D(A) = H^2(0,1;\R^2) \cap H^1_0(0,1 ; \R^2).
\end{gather*}
The following proposition is proved in~\cite[Proposition 2.1, Lemma 2.3]{AKBGBdT_JMAA16}.
\begin{prop} \label{Prop:Spectre_Syst_en_q}
The spectrum of $-A^*$ is given by the sequence of eigenvalues $\lambda_k = k^2 \pi^2$. 
Let
\[
\varphi_{k,1} = \begin{pmatrix} 0 \\ \varphi_k \end{pmatrix}, 
\qquad
\varphi_{k,2} = \begin{pmatrix} \varphi_k \\ \psi_k \end{pmatrix},
\]
with $(\lambda_k, \varphi_k)$ is given by~\eqref{EltsPropresClassiques} and $\psi_k$ is the unique solution of 
\begin{equation}
\left\{
\begin{aligned}
& -\frac{\md^2}{\md x^2} \psi(x) - \lambda_k \psi(x) = \left( I_k(q) -q(x) \right) \varphi_k(x), \quad &x \in (0,1),
\\
&\psi(0) = \psi(1) = 0,
\\
&\int_0^1 \psi(x) \varphi_k(x) \md x = 0.
\end{aligned}
\right.
\end{equation}
Then, 
\[
(A^* + \lambda_k) \varphi_{k,1} = 0, \qquad 
(A^* + \lambda_k) \varphi_{k,2} = -I_k(q) \varphi_{k,1},
\]
and the family $\{ \varphi_{k,1} , \varphi_{k,2} : k \in \N^* \}$ forms a Riesz basis of $H$.
\end{prop}

Moreover it is proved in~\cite[Theorem 1.3]{AKBGBdT_JMAA16} that if
\begin{equation} \label{Condition_Ik(q)_I1k(q)}
|I_k(q)| + |I_{1,k}(q)| \neq 0, \quad \forall k \in \N^*,
\end{equation}
which is a necessary and sufficient condition for approximate controllability, then there is a minimal time for null controllability of system~\eqref{Syst_en_q} in $H$ given by
\[
T_{min} = \limsup\limits_{k \to +\infty} \: \frac{\min \{ - \ln|I_k(q)|, -\ln|I_{1,k}(q)| \}}{\lambda_k}.
\]
Here, it is implicitely assumed that $\ln(0)=-\infty$.

\medskip
From Proposition~\ref{Prop:Spectre_Syst_en_q} it comes that if we assume that $I_k(q )=0$ for any $k \in \N^*$ then the eigenvalue $\lambda_k$ has geometric multiplicity equal to $2$. Thus, in this case, system~\eqref{Syst_en_q} fits into the setting of Theorem~\ref{ThCondensationNulle} (recall that $\operatorname{c}(\Lambda)=0$) and we deduce that
\[
T^* = T_{min} = \limsup\limits_{k \to +\infty} \: \frac{-\ln|I_{1,k}(q)|}{\lambda_k}.
\]

\begin{rmq}
In this example, the general case where $I_k(q) \neq 0$ for some $k \in \N^*$, does not fit into our abstract setting. One needs to take into account generalised eigenfunctions to obtain a complete sequence but assumption~\eqref{HypoStructureBJordan} is not satisfied by the internal control operator.
Nevertheless, this general case is discussed in Sec.~\ref{Subsec:FurtherSyst_q} where we still prove that $T^* = T_{min}$ (see Proposition~\ref{Prop:Syst_en_q}).
\end{rmq}

\subsubsection{A parabolic cascade system with a variable coupling coefficient: boundary control}
\label{ExempleSyst_q_bord}

Let $q \in L^\infty(0,1;\R)$. Let us consider the following system studied in~\cite{AKBGBdT_JMAA16}
\begin{equation} \label{Syst_en_q_bord}
\left\{
\begin{aligned}
& \partial_t y + \begin{pmatrix}  -\partial_{xx} & 0 \\ 0 & -\partial_{xx} \end{pmatrix} y + \begin{pmatrix} 0 & q(x) \\ 0 & 0 \end{pmatrix} y  = 0, \quad &(t,x) \in (0,T)\times(0,1),
\\
& y(t,0) = \begin{pmatrix} 0 \\ u(t) \end{pmatrix}, \quad y(t,1) = 0,  &t \in (0,T),
\\
&y(0,x) = y_0(x),    &x \in (0,1),
\end{aligned}
\right.
\end{equation}
We keep the notations introduced in the previous example. It is proved in~\cite[Theorem 1.1]{AKBGBdT_JMAA16} that if
\begin{equation} \label{Condition_Ik(q)}
|I_k(q)|  \neq 0, \quad \forall k \in \N^*,
\end{equation}
which is a necessary and sufficient condition for approximate controllability, then there is a minimal time for null controllability of system~\eqref{Syst_en_q_bord} in $H^{-1}(0,1 ; \R^2)$ with controls in $L^2(0,T;\R)$ given by
\[
T_{min} = \limsup\limits_{k \to +\infty} \: \frac{- \ln|I_k(q)|}{\lambda_k}.
\]
From Proposition~\ref{Prop:Spectre_Syst_en_q} it comes that system~\eqref{Syst_en_q_bord} fits into the setting of Proposition~\ref{PropJordan} (the assumption~\eqref{HypoStructureBJordan} being automatically satisfied as we are dealing with a scalar control). For a precise definition of the operator $B$ as well as its admissibility we refer for example to~\cite[Sec. 6.2]{AKBGBDT_condensation}. Thus $T_{min} \in [T^*, 2 T^*]$.

Moreover it comes that 
\[
\| B^* \varphi_{k,1} \|_{U} = |\varphi_k'(0)| = \sqrt{2} k \pi.
\]
Thus, from Corollary~\ref{CoroJordan}, one has $T^* = T_{min}$.

%%%%%%%%%%%%%%%%%%%%%%%%%%%%%%%%%%%%%%%%%%%%%%%%%%%%%%%%%%%%%%%%%%%%%%%%%%%%%%%%%%%%%%%%
\section{Condensation of the spectrum}
\label{Sec:PreuveCondensation}
%%%%%%%%%%%%%%%%%%%%%%%%%%%%%%%%%%%%%%%%%%%%%%%%%%%%%%%%%%%%%%%%%%%%%%%%%%%%%%%%%%%%%%%%

In this section we prove Theorems~\ref{ThCondensation} and~\ref{ThCondensationGeneral}. The sketch of proof is the same as Theorem~\ref{ThCondensationNulle} except that the condensation of eigenvalues has to be taken into account while estimating the biorthogonal family.

%---------------------------------------------------------------------------------------
\subsection{Condensation of eigenvalues given by the Bohr index}
%---------------------------------------------------------------------------------------

The case $T<T^*$ follows  directly from Theorem~\ref{ThDuyckaertsMiller}. 

Assume that $T^* < +\infty$. 
We start by proving that in this setting the condensation index is smaller than $T^*$ (and thus finite).
\begin{prop} \label{PropIndiceCondensationFini}
Assume that the hypotheses of Theorem~\ref{ThCondensation} hold. Then, if $T^* < +\infty$ it comes that
\[
\mathrm{c}(\Lambda) \leq T^*.
\]
\end{prop}

\begin{proof}
From Hypothesis~\eqref{HypoStructureB} it comes that there exist $\alpha_k, \alpha_j \in \C$ such that 
\[
\alpha_k \varphi_k + \alpha_j \varphi_j \in \mathrm{Ker}(B^*).
\]
Inequality~\eqref{InegResolvante} implies that $B^* \varphi_k \neq 0$. Thus, there exists $\alpha_{k,j} \in \C$ such that
\[
v_{k,j} := \alpha_{k,j} \varphi_k + \varphi_{j} \in \mathrm{Ker} (B^*).
\]
Applying \eqref{InegResolvante} it comes that for any $\varepsilon >0$ there exists $C_\varepsilon > 0$ such that for any $k \in \N^*$ and any $j \in \N^*$,
\begin{equation} \label{BohrTemp}
\| v_{k,j} \|^2 \leq C_\varepsilon e^{2 (T^* +\varepsilon) \Re(\lambda_k)} \frac{\|(A^*+\lambda_k) v_{k,j} \|^2}{\Re(\lambda_k)^2}.
\end{equation}
A direct computation leads to $(A^* + \lambda_k) v_{k,j} = (\lambda_{k}-\lambda_j) \varphi_{j}$. Thus, as we assumed that the eigenfunctions form a Riesz basis, there exists $c >0$ such that for any $k, j \in \N^*$
\begin{equation} \label{Norme_VP_orthogonalite}
\| v_{k,j} \|^2 \geq c (1 + \alpha_{k,j}^2) \geq c.
\end{equation}
This leads to
\begin{equation} \label{EstimeeIndiceBohr}
|\lambda_{k} - \lambda_j| \geq C_\varepsilon \Re(\lambda_k) e^{- \Re(\lambda_k) (T^*+\varepsilon)}.
\end{equation}
Thus, estimates~\eqref{BohrTemp},~\eqref{Norme_VP_orthogonalite} and~\eqref{EstimeeIndiceBohr} imply $\mathrm{Bohr}(\Lambda) \leq T^*$. The conclusion follows from the assumption $\mathrm{Bohr(\Lambda)}=\mathrm{c}(\Lambda)$. 

\end{proof}

\begin{rmq}
As seen in the previous proof, inequality~\eqref{InegResolvante} allows to estimate, under the structural assumption~\eqref{HypoStructureB}, the Bohr index of the sequence of eigenvalues. However, it is not clear if inequality~\eqref{InegResolvante} leads to an estimate for the condensation index in the general case. This is the main reason why we imposed the extra condition $\mathrm{Bohr}(\Lambda)= \mathrm{c}(\Lambda)$. 
\end{rmq}

One should notice that the structural assumption~\eqref{HypoStructureB} is necessary for Proposition~\ref{PropIndiceCondensationFini} to hold. Consider the following system
\begin{equation*}
\left\{
\begin{aligned}
&\partial_ t y = \begin{pmatrix} \partial_{xx} & 0 \\ 0 & d \partial_{xx} \end{pmatrix} y, \quad (t,x) \in (0,T) \times (0,1),
\\
&y(t,0) = \begin{pmatrix} u_1(t) \\ u_2(t) \end{pmatrix}, \quad y(t,1)=0,
\\
&y(0,x) = y_0(x),
\end{aligned}
\right.
\end{equation*}
with $d>0$.
It is direct to notice that controllability in $H^{-1}(0,1;\R^2)$ holds in any time $T>0$ and thus $T^*=0$. However, from~\cite[Proposition 6.20]{AKBGBDT_condensation}, one can choose $d$ to prescribe any value in $[0,+\infty]$ for the condensation index of the spectrum of the underlying operator.

%\begin{rmq} \label{RmqT0=condensation} 
%As stated in Remark~\ref{RmqT0=condensationAnnonce}, if we moreover consider a scalar control, then it comes that $T^* = c(\Lambda)$. Indeed from~\cite{AKBGBDT_condensation}, it comes that in this setting if $T> c(\Lambda)$ then system~\eqref{SystControlAbstrait} is null controllable in time $T$. Then, Theorem~\ref{ThDuyckaertsMiller} implies that $T \geq T^*$ and thus $c(\Lambda) \geq T^*$. The reverse inequality is proved in Proposition~\ref{PropIndiceCondensationFini}.
%\end{rmq}

\begin{proof}[Proof of Theorem~\ref{ThCondensation}.]
Let $T>c(\Lambda)$. 
As we assumed that the eigenvalues are all simple the strategy is exactly the one detailed in Sec.~\ref{Subsec:Strategy}. Indeed the control given by~\eqref{SolutionPbMomentAbstrait} is a formal solution of the moment problem~\eqref{PbMomentAbstrait2}. Thus, in view of~\eqref{SolutionPbMomentAbstrait}, the only remaining point is to prove that the following series 
\[
\sum_{k \in \N^*} e^{-\lambda_k T} \frac{\langle y_0, \varphi_k \rangle_{{ }}}{\| B^* \varphi_k \|_{U}^2} q_k(T-t) B^* \varphi_k
\]
converges in $L^2(0,T;U)$.

\medskip
Recall that, from Proposition~\ref{FamilleBiorthogonale}, we deduce that for any $\varepsilon >0$ there exists $C_\varepsilon >0$ such that for any $k \in \N^*$
\begin{equation} 
\| q_k \|_{L^2(0,T)} \leq C_\varepsilon e^{\Re(\lambda_k)(\mathrm{c}(\Lambda) + \varepsilon)}.
\end{equation}
%Thus, using Proposition~\ref{PropIndiceCondensationFini}, it comes that 
%\begin{equation} \label{EstimeeBiorthoCondensation}
%\| q_k \|_{L^2(0,T)} \leq C_\varepsilon e^{\Re(\lambda_k)(T^* + \varepsilon)}.
%\end{equation}
As a consequence, we obtain that
\begin{align*}
&\left\| e^{-\lambda_k T} \frac{\langle y_0, \varphi_k \rangle_{{ }}}{\| B^* \varphi_k \|_{U}^2} q_k(T-\cdot) B^* \varphi_k \right\|_{L^2(0,T;U)}
\\
&\leq C_\varepsilon  \frac{|\langle y_0, \varphi_k\rangle_{{ }}|}{e^{\varepsilon \Re(\lambda_k)} \|B^* \varphi_k\|_{U}} e^{-\Re(\lambda_k)(T-c(\Lambda)-2\varepsilon)}.
\end{align*}
As $T>c(\Lambda)$ the choice of $\varepsilon$ sufficiently small leads to null controllability in time $T$. 

Finally it comes that null controllability holds in time $T>c(\Lambda)$. From Theorem~\ref{ThDuyckaertsMiller} null controllability does not hold in time $T<T^*$. Recall that Proposition~\ref{PropIndiceCondensationFini} yields $c(\Lambda) \leq T^*$ and thus ends the proof.

\end{proof}

\begin{rmq} \label{RmqTminBohr}
Notice that if in the previous proof, instead of assuming~\eqref{ActionControleSuffisante}, we had applied inequality~\eqref{InegResolvante} with $v=\varphi_k$ and $\lambda=\lambda_k$ we would have obtained the following estimate
\[ 
\| B^* \varphi_k \|_{U} \geq C_\varepsilon \sqrt{\Re(\lambda_k)} e^{-\Re(\lambda_k) (T^*+\varepsilon)}.
\]
This would have led to controllability in time $T>T^*+c(\Lambda)$. This result is contained in Theorem~\ref{ThCondensationGeneral} proved in the next subsection.
\end{rmq}

\begin{rmq} \label{Rmq:HypoStructureB_precisee}
A careful inspection of the proof of Theorem~\ref{ThCondensation} shows that assumption~\eqref{HypoStructureB} can be replaced by the less explicit assumption
\begin{equation} \label{HypoStructureB_precisee}
\mathrm{Ker}(B^*) \cap \mathrm{Span}(\varphi_{k_n}, \varphi_{j_n}) \neq \{0\}, \quad \forall n \in \N^*,
\end{equation}
where the sequence $(k_n)_{n \in \N^*}$ is such that
\[
\lim\limits_{n \to +\infty} \: \frac{-\ln \inf\limits_{j \neq k_n }|\lambda_{k_n}-\lambda_j|}{\Re(\lambda_{k_n})}=
\limsup\limits_{k \to +\infty} \: \frac{-\ln \inf\limits_{j \neq k }|\lambda_{k}-\lambda_j|}{\Re(\lambda_k)}
\]
and the sequence $(j_n)_{n \in \N^*}$ is such that
\[
|\lambda_{k_n} - \lambda_{j_n}| = \inf\limits_{j \neq k_n} |\lambda_{k_n} - \lambda_j|.
\]
\end{rmq}

%---------------------------------------------------------------------------------------
\subsection{A qualitative result in the general setting}
%---------------------------------------------------------------------------------------

We prove in this subsection the qualitative result in the most general setting considered, namely Theorem~\ref{ThCondensationGeneral}. The proof follows the lines of the proof of Theorem~\ref{ThCondensationNulle} except that we take into account the effect of the condensation index in the estimate of the biorthogonal family (in time). 

\begin{proof}[Proof of Theorem~\ref{ThCondensationGeneral}.]

The lack of null controllability when $T<T^*$ follows directly from Theorem~\ref{ThDuyckaertsMiller}.
We thus focus on the case $T^* <+\infty$ and we set $T>T^*+\operatorname{c}(\Lambda)$.

$\bullet$ To clarify the presentation let us start with the simplified case where the eigenvalues are geometrically simple, i.e. $r_k=1$.

As already used the control given by~\eqref{SolutionPbMomentAbstrait} is a formal solution of the moment problem~\eqref{PbMomentAbstrait2}. Thus, in view of~\eqref{SolutionPbMomentAbstrait}, the only thing to prove is that the following series 
\[
\sum_{k \in \N^*} e^{-\lambda_k T} \frac{\langle y_0, \varphi_k \rangle_{{ }}}{\| B^* \varphi_k \|_{U}^2} q_k(T-t) B^* \varphi_k
\]
converges in $L^2(0,T;U)$. 

For any $\varepsilon>0$, inequality~\eqref{InegResolvante}, implies that 
\[
\| B^* \varphi_k \|_{U} \geq C_\varepsilon \sqrt{\Re(\lambda_k)} e^{-\Re(\lambda_k)(T^*+\varepsilon)}
\]
and Proposition~\ref{FamilleBiorthogonale} implies that
\[
\| q_k(T-\cdot) \|_{L^2(0,T)} \leq C_\varepsilon e^{\Re(\lambda_k) (\mathrm{c}(\Lambda) + \varepsilon)}.
\]
Thus, the formal solution is actually a solution of the moment problem in $L^2(0,T;U)$ as soon as $T > T^* + \mathrm{c}(\Lambda)$.

$\bullet$ We now turn to the case of geometrically multiple eigenvalues. 
As in the proof of Theorem~\ref{ThCondensationNulle} we define the sets
\[
\Sigma_S = \{ k \in \N^* ; r_k=1\}, \qquad \Sigma_M = \{ k \in \N^* ; r_k \geq 2\}
\]
and the control $u = u_1 + u_2$ where $u_1$ is given by~\eqref{FormeControleVPSimple} and 
\[
u_2(t) = - \sum_{k \in \Sigma_M} q_k(T-t) \left( \sum_{1 \leq i \leq r_k} \langle y_0, \varphi_{k,i} \rangle_{{ }} e^{-\lambda_k T} \Psi_{k,i} \right),
\]
with $(\Psi_{k,1},\dots, \Psi_{k,r_k})$ the biorthogonal family to $(B^* \varphi_{k,1}, \dots, B^* \varphi_{k,r_k})$ given by Lemma~\ref{LemmeBiorthoEspace}. The same computations as in the proof of Theorem~\ref{ThCondensationNulle} imply the same estimate~\eqref{EstimeeFamilleBiorthoEspace} for this biorthogonal family. 
This, together with Proposition~\ref{FamilleBiorthogonale}, directly gives that $u \in L^2(0,T;U)$ if $T> T^* + \mathrm{c}(\Lambda)$ and $u$ is a solution of the moment problem.

\end{proof}

%---------------------------------------------------------------------------------------
\subsection{Examples} \label{Subsec:ExamplesCondensation}
%---------------------------------------------------------------------------------------

\subsubsection{An academic example}

Let $(\lambda_k,\varphi_k)$ be defined as in~\eqref{EltsPropresClassiques}. Let $f : \R \to \R$ satisfying
\[
0 < |f(s)| < s, \: \forall s \quad \text{and} \quad f(s) = \underset{s \to +\infty}{o}(s).
\]
Consider the operator $L_f$ defined by
\[
L_f = \sum_{k \in \N^*} f(\lambda_k) \langle \cdot , \varphi_k \rangle \varphi_k.
\]
Let us consider the following system studied in~\cite[Sec.~2.2]{AKBDK05}
\begin{equation} \label{SystIlya}
\left\{
\begin{aligned}
& \partial_t y = \begin{pmatrix} \partial_{xx} & L_f \\ L_f & \partial_{xx} \end{pmatrix} y + \begin{pmatrix}
0 \\ u(t) \end{pmatrix}  \quad & (t,x) \in (0,T) \times (0,1),
\\
& y(t,0) = y(t,1) = 0,       & t\in (0,T),
\\
&y(0,x) = y_0(x),   &x \in (0,1),
\end{aligned}
\right.
\end{equation}

Let $H=L^2(0,1 ; \R^2)$ and $U= \R$. Define the operator $A$ by 
\begin{gather*}
A y = \begin{pmatrix} \partial_{xx} & L_f \\ L_f & \partial_{xx} \end{pmatrix} y,
\\
D(A) = H^2(0,1 ; \R^2) \cap H^1_0(0,1 ; \R^2),
\end{gather*}
and the bounded control operator $B = \begin{pmatrix} 0 \\ 1\end{pmatrix}$.

The spectrum of the operator $-A^*$ is given by $\Lambda =\{ \lambda_k \pm f(\lambda_k) ; k \in \N^* \}$ and the associated eigenfunctions are given by 
\[
\varphi_k^\pm = \frac{1}{\sqrt{2}} \begin{pmatrix} \varphi_k \\ \mp \varphi_k \end{pmatrix}.
\]
Moreover from~\cite{AKBDK05} it can be proved that there exists a minimal null control time for system~\eqref{SystIlya} given by
\[
T_{min} = \limsup\limits_{k \to +\infty} \: \frac{-\ln|f(\lambda_k)|}{\lambda_k}.
\]

For $k$ sufficiently large the eigenvalues are ordered in the following order
\[
\lambda_k -f(\lambda_k) < \lambda_k + f(\lambda_k) < \lambda_{k+1} - f(\lambda_{k+1}).
\]
Thus the spectrum $\Lambda$ can be split into the following condensation grouping $G_k =\{ \lambda_k \pm f(\lambda_k) \}$. 
From~\cite{Shackell69} we deduce that
\[
\mathrm{c}(\Lambda) = \mathrm{Bohr}(\Lambda) = \limsup\limits_{k \to +\infty} \: \frac{-\ln|f(\lambda_k)|}{\lambda_k}.
\]

Given the particular structure of the eigenfunctions it comes that 
\[
\varphi_k^+ + \varphi_k^- = \begin{pmatrix} \sqrt{2} \varphi_k \\ 0 \end{pmatrix}
\]
and thus $B^*(\varphi_k^+ + \varphi_k^-) = 0$. According to Remark~\ref{Rmq:HypoStructureB_precisee} it comes that system~\eqref{SystIlya} fits into the setting of Theorem~\ref{ThCondensation}. Thus, it comes that
\[
T^* = T_{min} = \limsup\limits_{k \to +\infty} \: \frac{-\ln|f(\lambda_k)|}{\lambda_k}.
\]

\begin{rmq}
This example is generalised in Sec.~\ref{Subsec:IlyaGeneral} in any space dimension. Though this generalisation will not fit into our abstract settings it will exhibit a minimal null control time still given by $T^*$.
\end{rmq}

\subsubsection{A parabolic cascade system with different diffusions: boundary control}
\label{ExempleSystDeuxDiffusions}

Let us consider  the following system studied in~\cite[Sec. 6.2]{AKBGBDT_condensation}
\begin{equation} \label{SystDeuxDiffusions}
\left\{
\begin{aligned}
& \partial_t y = \begin{pmatrix} \partial_{xx} & 0 \\ 0 & d \partial_{xx} \end{pmatrix} y + \begin{pmatrix} 0 & 1 \\ 0 & 0 \end{pmatrix} y,  \quad & (t,x) \in (0,T) \times (0,1),
\\
& y(t,0) = \begin{pmatrix} 0 \\ u(t) \end{pmatrix}, \quad y(t,1) = 0,       & t\in (0,T),
\\
&y(0,x) = y_0(x),   &x \in (0,1),
\end{aligned}
\right.
\end{equation}
with $d \in (0,+\infty)$, $d\neq 1$.

To fit into our abstract setting, we set $H=H^{-1}(0,1; \R^2)$, $U=\R$ and 
\begin{gather*}
A = \begin{pmatrix} \frac{\md^2}{\md x^2} & 0 \\ 0 & d \frac{\md^2}{\md x^2} \end{pmatrix} + \begin{pmatrix} 0 & 1 \\ 0 & 0 \end{pmatrix},
\\
D(A) = H^2(0,1;\R^2) \cap H^1_0(0,1 ; \R^2).
\end{gather*}
The precise definition of the operator $B$ as well as its admissibility are given in~\cite[Sec. 6.2]{AKBGBDT_condensation}.

Let $(\lambda_k, \varphi_k)$  be as in~\eqref{EltsPropresClassiques}. The following proposition is proved in~\cite[Propositions 6.11 and 6.12]{AKBGBDT_condensation}. 
\begin{prop} \label{Prop:Spectre_Syst_en_d}
The spectrum of $-A^*$ is given by the sequence of eigenvalues $\Lambda= \{\lambda_{k,1}= \lambda_k, \lambda_{k,2} = d \lambda_k ; k \in \N^*\}$. 
These eigenvalues are simple if and only if $\sqrt{d} \not\in \mathbb{Q}$. 

If $\sqrt{d} \not\in \mathbb{Q}$, the associated eigenfunctions are given by
\[
\varphi_{k,1} = \sqrt{\lambda_k} \begin{pmatrix} \varphi_k \\ \psi_k \end{pmatrix}, 
\qquad
\varphi_{k,2} = \sqrt{\lambda_k} \begin{pmatrix} 0 \\ \varphi_k \end{pmatrix},
\]
with $\psi_k = \frac{1}{\lambda_k(d-1)} \varphi_k$.

The family $\{ \varphi_{k,1} , \varphi_{k,2} : k \in \N^* \}$ forms a Riesz basis of $H$ and 
\[
B^* \varphi_{k,1} = \frac{\sqrt{2}}{d-1}, \qquad B^* \varphi_{k,2} = \sqrt{2} \lambda_k.
\]
\end{prop}

Moreover it is proved in~\cite[Theorem 6.17]{AKBGBDT_condensation} that if
\begin{equation} \label{Condition_sqrt(d)}
\sqrt{d} \not\in \mathbb{Q}
\end{equation}
which is a necessary and sufficient condition for approximate controllability, then there is a minimal time for null controllability of system~\eqref{SystDeuxDiffusions} in $H$ given by
\[
T_{min} = \mathrm{c}(\Lambda).
\]

\medskip
From Proposition~\ref{Prop:Spectre_Syst_en_d}, it comes that the only remaining point to prove to fit system~\eqref{SystDeuxDiffusions} into the setting of Theorem~\ref{ThCondensation} is $\mathrm{c}(\Lambda)=\mathrm{Bohr}(\Lambda)$. 

This property is proved in~\cite{Samb17}. For the sake of completeness, we reproduce the proof in~\ref{AnnexeBohr}. 

Thus we deduce that
\[
T^* = T_{min} = \mathrm{c}(\Lambda).
\]

\subsubsection{A parabolic cascade system with different diffusions: pointwise control}
\label{ExSystDeuxDiffusionsControlePonctuel}
Let us consider  the following system studied in~\cite[Sec. 6.3]{AKBGBDT_condensation}
\begin{equation} \label{SystDeuxDiffusions_ControlePonctuel}
\left\{
\begin{aligned}
& \partial_t y = \begin{pmatrix} \partial_{xx} & 1 \\ 0 & d \partial_{xx} \end{pmatrix} y + \begin{pmatrix}
0 \\ \delta_{x_0} u(t) \end{pmatrix}  \quad & (t,x) \in (0,T) \times (0,1),
\\
& y(t,0) =  y(t,1) = 0,       & t\in (0,T),
\\
&y(0,x) = y_0(x),   &x \in (0,1),
\end{aligned}
\right.
\end{equation}
with $d \in (0,+\infty)$, $d\neq 1$.

The abstract setting for the operator $A$ is the same as the one introduced in Sec.~\ref{ExempleSystDeuxDiffusions}. Thus the spectral analysis is the one given in Proposition~\ref{Prop:Spectre_Syst_en_d}. The abstract setting for the operator $B$ is the one introduced in Sec.~\ref{ExempleDolecki}.

Moreover it is proved in~\cite[Theorem 6.31]{AKBGBDT_condensation} that if
\begin{equation} \label{Condition_sqrt(d)_Dolecki}
\sqrt{d} \not\in \mathbb{Q} \quad \text{ and } \quad \varphi_k(x_0) \neq 0, \: \forall k \in \N^*
\end{equation}
which is a necessary and sufficient condition for approximate controllability, then there is a minimal time for null controllability of system~\eqref{SystDeuxDiffusions_ControlePonctuel} in $H$ given by
\[
T_{min} = \max\limits_{i=1, 2} \: \limsup\limits_{k \to +\infty} \: \frac{-\ln|\varphi_k(x_0)| - \ln |E'(\lambda_{k,i})|}{\lambda_{k,i}}.
\]

\medskip
Thus, if $d$ is such that $c(\Lambda) <+\infty$, system~\eqref{SystDeuxDiffusions_ControlePonctuel} directly fits into the setting of Theorem~\ref{ThCondensationGeneral} and we deduce that 
\[
T^* \leq T_{min} \leq T^* + c(\Lambda). 
\]

%%%%%%%%%%%%%%%%%%%%%%%%%%%%%%%%%%%%%%%%%%%%%%%%%%%%%%%%%%%%%%%%%%%%%%%%%%%%%%%%%%%%%%%%
\section{Further results}
\label{Sec:AutresCas}
%%%%%%%%%%%%%%%%%%%%%%%%%%%%%%%%%%%%%%%%%%%%%%%%%%%%%%%%%%%%%%%%%%%%%%%%%%%%%%%%%%%%%%%%

This section is dedicated to further results. First we prove that for some examples of parabolic problems exhibiting a (positive finite) minimal null control time (that does not fit into the abstract setting developed in this article), this minimal time is still given by $T^*$. Then, we propose an academic example (inspired by~\eqref{SystIlya}) to give a multi-dimensional setting where the minimal time is given by $T^*$. Unlike the rest of this article, the proof does not rely on the moment method. 
Finally we end this article with a (counter-) example showing that hypothesis~\eqref{ComportementSpectre} is not a technical assumption due to the use of moment method but is deeply related to the validity of the presented results.

%---------------------------------------------------------------------------------------
\subsection{Known examples of minimal time not covered by the abstract setting.} 
%---------------------------------------------------------------------------------------

\subsubsection{A parabolic cascade system with a variable coupling coefficient: the general case of internal control}
\label{Subsec:FurtherSyst_q}

Let us consider the setting and notations detailed in Sec.~\ref{ExempleSyst_q_CasParticulier} for the study of system~\eqref{Syst_en_q} but with the general hypothesis
\[
|I_k(q)| + |I_{1,k}(q)| \neq 0.
\]
Recall that, in this setting, the eigenvalue $\lambda_k$ could have algebraic multiplicity equal to $2$. One has,
\begin{prop} \label{Prop:Syst_en_q}
Let $T^*$ be the minimal time defined by~\eqref{DefTempsMinimal}. Then 
\[
T^* = T_{min} = \limsup\limits_{k \to +\infty} \: \frac{\min \{ - \ln|I_k(q)|, -\ln|I_{1,k}(q)| \}}{\lambda_k}.
\]
\end{prop}

\begin{proof}
Assume that $T>T_{min}$. Then, from~\cite{AKBGBdT_JMAA16}, system~\eqref{Syst_en_q} is null controllable. Thus, Theorem~\ref{ThDuyckaertsMiller} implies that inequality~\eqref{InegResolvante} holds and thus $T \geq T^*$ which gives $T^* \leq T_{min}$.

Assume now that $T$ is such that inequality~\eqref{InegResolvante} holds. From the definition of $T_{min}$ we deduce that, for a subsequence still denoted $k$, for any $\varepsilon>0$, there exists $K \in \N^*$ such that 
\begin{equation} \label{Estimee_Ik(q)}
\max \{ |I_k(q)|, |I_{1,k}(q)| \} \leq e^{-\lambda_k(T_{min}-\varepsilon)}, \quad \forall k \geq K.
\end{equation}
We will use the following lemma which gives details on the structure of the generalised eigenfunctions. Its proof is given in~\cite[Proposition 2.6]{AKBGBdT_JMAA16}.
\begin{lemme}
For any $k \in \N^*$, there exists $\tau_k \in \R$ such that
\[
\psi_k(x) = \tau_k \varphi_k(x) + \xi_k(x), \qquad \forall x \in \omega,
\]
with $\psi_k$ given in Proposition~\ref{Prop:Spectre_Syst_en_q} and
\[
\|\xi_k\|_{L^2(\omega)} \leq C \left( |I_k(q)| + |I_{1,k}(q)| \right).
\]
\end{lemme}
Then, we apply inequality~\eqref{InegResolvante} with $\lambda=\lambda_k$ and $v=\varphi_{k,2} - \tau_k \varphi_{k,1} = \begin{pmatrix} \varphi_k \\ \psi_k - \tau_k \varphi_k \end{pmatrix}$.
Using the previous lemma and Proposition~\ref{Prop:Spectre_Syst_en_q}, we deduce that,
\begin{gather*}
\|v\|^2 = 1 + \|\ \psi_k - \tau_k \varphi_k \|^2 \geq 1,
\\
\| (A^*+\lambda) v\|^2 = \| I_k(q) \varphi_{k,1} \|^2 = |I_k(q)|^2,
\\
\| B^* v \|_{U}^2 = \| \psi_k - \tau_k \varphi_k \|_{L^2(\omega)} = \| \xi_k \|_{L^2(\omega)}^2 \leq C  \left( |I_k(q)|^2 + |I_{1,k}(q)|^2 \right),
\end{gather*}
it comes that 
\[
1 \leq \frac{C}{\lambda_k} e^{2 \lambda_k T}  \left( |I_k(q)|^2 + |I_{1,k}(q)|^2 \right).
\]
Using the estimate~\eqref{Estimee_Ik(q)} it comes that
\[
1 \leq \frac{C}{\lambda_k} e^{2 \lambda_k (T-(T_{min} -\varepsilon))}.
\]
Letting $k$ go to $+\infty$ implies $T \geq T_{min} - \varepsilon$ and thus $T^* \geq T_{min}$ which ends the proof.

\end{proof}

\subsubsection{Extension to cascade systems with a first order coupling term}

In~\cite{Duprez_Tmin2017} the results of~\cite{AKBGBdT_JMAA16} concerning the controllability of system~\eqref{Syst_en_q} and system~\eqref{Syst_en_q_bord} are extended (using the same technics) to the following systems
\begin{equation*} 
\left\{
\begin{aligned}
& \partial_t y + \begin{pmatrix}  -\partial_{xx} & 0 \\ 0 & -\partial_{xx} \end{pmatrix} y + \begin{pmatrix} 0 & p(x) \partial_x + q(x) \\ 0 & 0 \end{pmatrix} y  = \begin{pmatrix} 0 \\ u \mathbf{1}_\omega \end{pmatrix},  
\\
& y(t,0) = y(t,1) = \begin{pmatrix} 0 \\ 0 \end{pmatrix},  
\\
&y(0,x) = y_0(x),
\end{aligned}
\right.
\end{equation*}
and
\begin{equation*}
\left\{
\begin{aligned}
& \partial_t y + \begin{pmatrix}  -\partial_{xx} & 0 \\ 0 & -\partial_{xx} \end{pmatrix} y + \begin{pmatrix} 0 & p(x) \partial_x + q(x) \\ 0 & 0 \end{pmatrix} y  = 0,
\\
& y(t,0) = \begin{pmatrix} 0 \\ u(t)\end{pmatrix},
\\
& y(t,1) = 0,   
\\
&y(0,x) = y_0(x).
\end{aligned}
\right.
\end{equation*}
The results stated in Sec.~\ref{ExempleSyst_q_CasParticulier}, Sec.~\ref{ExempleSyst_q_bord} and Proposition~\ref{Prop:Syst_en_q} directly extend to these systems.

\subsubsection{Degenerate parabolic equation of Grushin-type}
\label{Subsec:FurtherGrushin}

In~\cite{BeauchardMillerMorancey_Tmin} it is proved that the following system
\begin{equation} \label{SystGrushin}
\left\{
\begin{aligned}
& \partial_t y - \partial_{x_1 x_1} y - x_1^2 \partial_{x_2 x_2} y = \mathbf{1}_\omega u(t,x_1,x_1), \quad &(t,(x_1,x_2)) \in (0,T) \times \Omega,
\\
&y(t,x_1,x_2) = 0, \quad &(t,(x_1,x_2)) \in (0,T) \times \partial \Omega,
\\
&y(0,x_1,x_2) = y_0(x_1,x_2),
\end{aligned}
\right.
\end{equation}
with $\Omega = (-1,1) \times (0,1)$ and $\omega=\Big( (-b,-a) \cup (a,b) \Big) \times (0,1)$ is null controllable in $L^2(\Omega)$ if and only if $T> T_{min}=\frac{a^2}{2}$. 
Contrarily to previously mentioned one-dimensional examples, this is a two dimensional degenerate parabolic problem.

Let us specify the abstract setting to express system~\eqref{SystGrushin} as~\eqref{SystControlAbstrait}. We set $H=L^2(\Omega)$ and $U=L^2(\omega)$.
For any $y, z \in C^\infty_0(\Omega)$, let
\[
(y,z) = \int_\Omega \left( \partial_{x_1} y  \partial_{x_1} z + x_1^2 \partial_{x_2} y  \partial_{x_2} \right) \md x_1 \md x_2
\]
$\| \cdot \|_V^2 = ( \cdot, \cdot )$ and $V = \overline{C^\infty_0(\Omega)}^{\|\cdot\|_V}$. The operator $A$ is defined by
\begin{gather*}
D(A) = \left\{ y \in V ; \exists c >0 \text{ such that } |(y,z)| \leq c \| z \|_{L^2(\Omega)}, \forall z \in V \right\},
\\
\langle Ay ,z \rangle = -(y,z), \quad \forall z \in V.
\end{gather*}
The operator $B$ is defined by $B = \mathbf{1}_\omega$.

Although this system does not fit into the considered settings (assumption~\eqref{ComportementSpectre} is not proved) the time $T^*$ is the minimal null control time.
\begin{prop} \label{Prop:Grushin}
Let $T^*$ be the minimal time defined by~\eqref{DefTempsMinimal} with the operators $A$ and $B$ defined above. Then,
\[
T^*=T_{min}=\frac{a^2}{2}.
\]
\end{prop}

\begin{rmq}
The existence of a minimal null control time for this system was already proved in~\cite{BeauchardCannarsaGuglielmi} in the case where $\omega =(a,b) \times (0,1)$ but its value is not yet known. The following proof requires the value of the minimal time to be known to to be exactly $\frac{a^2}{2}$. This is why we considered a control domain made of two symmetric strips $\omega=\Big( (-b,-a) \cup (a,b) \Big) \times (0,1)$.
\end{rmq}

The proof of Proposition~\ref{Prop:Grushin} given below follows closely the arguments of~\cite{BeauchardCannarsaGuglielmi}. The main motivation to study this example in this setting is to enlighten that the controllability criterion given by the definition of $T^*$~\eqref{DefTempsMinimal} goes beyond the framework studied in Theorems~\ref{ThCondensationNulle}, \ref{ThCondensation} and~\ref{ThCondensationGeneral}. 
\begin{proof}
Assume that $T> \frac{a^2}{2}$. Then as problem~\eqref{SystGrushin} is null controllable, Theorem~\ref{ThDuyckaertsMiller} implies that~\eqref{InegResolvante} holds and thus $T \geq T^*$ which gives $T^* \leq \frac{a^2}{2}$.

Assume now that $T$ is such that inequality~\eqref{InegResolvante} holds.
The following properties are proved in~\cite[Lemma 2, Lemma 4]{BeauchardCannarsaGuglielmi}. 

Let $\lambda_n$ be the smallest eigenvalue of the operator $-\frac{\md^2}{\md x^2} + (n \pi)^2 x^2$ on the domain $H^2(-1,1;\R) \cap H^1_0(-1,1;\R)$ and $v_n$ the associated $L^2(-1,1;\R)$ normalised eigenfunction. Then, there exists $c>0$ such that for any $n \in \N^*$, 
\[
n \pi \leq \lambda_n \leq n \pi +c.
\] 
Moreover, $v_n$ is even and 
\[
\int_a^b v_n^2(x) \md x \underset{n \to + \infty}{\sim}  \frac{e^{-a^2 n \pi}}{2a \pi \sqrt{n}}.
\]
Thus applying inequality~\eqref{InegResolvante} with $\lambda=\lambda_n$ and $v(x_1,x_2) = \sqrt{2} v_n(x_1) \sin(n \pi x_2)$ it comes that
\begin{align}
1 \leq C_T \frac{e^{2 \lambda_n T}}{\lambda_n} \| B^* v \|_{U}^2 &= \frac{2 C_T e^{2\lambda_n T}}{\lambda_n} \int_a^b v_n^2(x) \md x
\notag
\\
&\leq \frac{C}{\lambda_n} e^{2 n \pi T} \int_a^b v_n^2(x) \md x
\notag
\\
&\underset{n \to + \infty}{\sim} \frac{C}{\sqrt{n} \lambda_n} e^{2n \pi \left(T-\frac{a^2}{2}\right)}.
\label{InegGrushin}
\end{align}
Thus $T \geq \frac{a^2}{2}$. This leads to $T^* \geq \frac{a^2}{2}$ and ends the proof.

\end{proof}

\begin{rmq}
As it may not appear in the previous proposition (due to the formulation using the infimum), in this example, inequality~\eqref{InegResolvante} even captures the controllability property when the time is equal to the minimal null control time. Indeed inequality~\eqref{InegGrushin} implies that~\eqref{InegResolvante} cannot hold in time $T= \frac{a^2}{2}$. From Theorem~\ref{ThDuyckaertsMiller} it implies that system~\eqref{SystGrushin} cannot be controllable in time $T=\frac{a^2}{2}$.
\end{rmq}

%---------------------------------------------------------------------------------------
\subsection{An academic example in any dimension}
\label{Subsec:IlyaGeneral}
%---------------------------------------------------------------------------------------

In this section we give an abstract system (inspired by system~\eqref{SystIlya}), in any space dimension, which exhibits a minimal null control time given by $T^*$ defined by~\eqref{DefTempsMinimal}.

Assume that $H = X \times X$ with $X$ being a Hilbert space and that $U=X$. Assume that the operator $-A$ is positive self-adjoint on $H$ and denote by $\{\lambda_{k,i} ; i=1,2, k\in \N^*\}$ its sequence of eigenvalues sorted in the following order
\[
\cdots < \lambda_{k,1} < \lambda_{k,2} < \lambda_{k+1,1} < \cdots
\]
Assume moreover that the corresponding eigenfunctions have the following form
\[
\varphi_{k,i} = \varphi_k e_i,
\]
where $(e_1,e_2)$ is the canonical basis of $\R^2$ and $\{ \varphi_k ; k \in \N^*\}$ is a Hilbert basis of $X$. 

We consider the bounded control operator 
$B = \begin{pmatrix}  b_1 \\ b_2 \end{pmatrix} \in \R^2$.

\begin{prop} \label{Prop:AcademiqueGeneralise}
Assume that $A$ and $B$ are as defined above. Let $T^*$ be defined by~\eqref{DefTempsMinimal}.
\begin{itemize}
\item[$\bullet$] if $T^* >0$ and $T<T^*$ then system~\eqref{SystControlAbstrait} is not null controllable in time $T$;
\item[$\bullet$] if $T^* <+\infty$ and $T>T^*$ then system~\eqref{SystControlAbstrait} is null controllable in time $T$.
\end{itemize}
\end{prop}

\begin{proof}
Assume that $T^*<T<+\infty$. Inequality~\eqref{InegResolvante} with $v=\varphi_{k,i}$ and $\lambda= \lambda_{k,i}$ implies
$b_1 b_2 \neq 0$.

To prove null controllability of system~\eqref{SystControlAbstrait} we prove that there is $C>0$ such that for any $y_k^0 \in \R^2$ there exists $u_k \in L^2(0,T;\R)$ with $\|u_k\|_{L^2} \leq C \|y_k^0\|$ such that the associated solution of  
\begin{equation} \label{SysDimFinie}
\left\{
\begin{aligned}
&y_k'(t)=A_ky_k(t)+Bu_k(t),\quad \forall t>0\\
&y_k(0)=y_k^0, 
\end{aligned}
\right.
\end{equation}
with $A_k:=\begin{pmatrix} -\lambda_{k,1}&0\\0& -\lambda_{k,2}\end{pmatrix}$ satisfies $y_k(T)=0$. The Hilbert basis property of $\{ \varphi_k ; k \in \N^*\}$ will then end the proof. 

As it is classical (see for instance~\cite{Zabczyk}) this finite dimensional system is controllable and the minimal $L^2(0,T;\R)$ norm control is given by
\begin{equation} \label{cdimfinie}
u_k (t) =- B^*e^{A_k(T-t)}Q_{k,T}^{-1}e^{A_kT} y_k^0.
\end{equation}
where
\[
Q_{k,T}=\int_0^T e^{A_ks} BB^* e^{ A_ks} \md s.
\]
We now need to estimate the norm of this control. As $BB^*=\begin{pmatrix} b_1^2&b_1b_2\\b_1b_2&b_2^2\end{pmatrix}$, it comes that
\[
Q_{k,T}= T \begin{pmatrix} b_1^2\eta(-2T\lambda_{k,1})&b_1b_1\eta(-T(\lambda_{k,1}+\lambda_{k,2}))\\b_1b_1\eta(-T(\lambda_{k,1}+\lambda_{k,2})) & b_2^2\eta(-2T\lambda_{k,2}) \end{pmatrix},
\]
with  $\eta : s \mapsto \frac {e^{s}-1}{s}$. As,
\begin{equation}\label{normeminimale}
\Vert u_k\Vert_{L^2(0,T)}^2=\langle Q_{k,T}^{-1}e^{A_kT}y_k^0,e^{A_kT}y_k^0\rangle\leq \underset{1\leq i\leq 2}\max\,e^{-2\lambda_{i,k}T}\,\Vert Q_{k,T}^{-1}\Vert\Vert y_k^0\Vert^2,
\end{equation}
we will estimate the norm of $Q_{k,T}^{-1}$. Let $\sigma_k$ be its smallest eigenvalue. The following computations are closely related to those of~\cite[Sec.~2.2]{AKBDK05}.

As $Q_{k,T} \in \mathcal{M}_2(\R)$ we obtain
\begin{equation}\label{sigmadcc1161}     
\frac {2\det Q_{k,T}}{\operatorname{tr} \,Q_{k,T}}\geq   \sigma_k\geq   \frac {\det Q_{k,T}}{\operatorname{tr} \,Q_{k,T}},\quad \forall k\geq 1,
\end{equation}
where $\operatorname{tr}$ denotes the trace of a matrix.

For any $s_1, s_2 \in (-\infty,0)$ there exists $\bar s \in (s_1,s_2)$ such that   
\[
\eta(2s_1)\eta(2s_2)=\eta(s_1+s_2)^2e^{(s_2-s_1)^2  (\ln(\eta (\bar s)))''}.
\]
Explicit computations lead to 
\[
(\ln(\eta(s))''=\frac 1{s^2}-\frac{1}{(e^{\frac{s}{2}}-e^{-\frac{s}{2}})^2}
\]
and thus 
\[
\frac 1{s^2+12}\leq (\ln(\eta(s))''\leq \frac 1{s^2}.
\]
Let $s_1=-T\lambda_{k,1}$ and $s_2=-T\lambda_{k,2}$. Thus,
\[
\det Q_{k,T}= T^2 b_1^2b_2^2\eta(-2T\lambda_{k,1})\eta (-2T\lambda_{k,2})\Big(1-e^{-(\lambda_{k,1}-\lambda_{k,2})^2)\psi_k}\Big)
\]
with $\psi_k$ satisfying
\[
\frac 1{T^2\lambda_{k,2}^2+12}\leq \psi_k\leq \frac 1{4T^2\lambda_{k,1}^2}.
\]
Now the trace of $Q_{k,T}$ is given by
\[
\operatorname{tr} (Q_{k,T})= Tb_1^2\eta(-2T\lambda_{k,1})+Tb_2^2\eta(-2T\lambda_{k,2}).
\]
Gathering these informations leads to
\begin{align*}
\sigma_k &\geq \frac{Tb_1^2b_2^2\eta(-2T\lambda_{k,1})\eta (-2T\lambda_{k,2})\Big(1-\exp\left( \frac{-(\lambda_{k,1}-\lambda_{k,2})^2)}{  4\lambda_{k,1}^2 } \right) \Big)}{b_1^2\eta(-2T\lambda_{k,1})+b_2^2\eta(-2T\lambda_{k,2})}
\\
&\geq C \frac{Tb_1^2b_2^2\eta(-2T\lambda_{k,1})\eta (-2T\lambda_{k,2})\left( \lambda_{k,1}-\lambda_{k,2} \right)^2 }{\lambda_{k,1}^2 \left(b_1^2\eta(-2T\lambda_{k,1})+b_2^2\eta(-2T\lambda_{k,2}) \right)}
\end{align*}
We now use inequality~\eqref{InegResolvante} to estimate $|\lambda_{k,2}-\lambda_{k,1}|$ in the same spirit as in the proof of Proposition~\ref{PropIndiceCondensationFini}. Applying inequality~\eqref{InegResolvante} with $v = \varphi_{k,2} - \frac{b_2}{b_1} \varphi_{k,1}$ and $\lambda=\lambda_{k,1}$ we obtain for any $\varepsilon>0$,
\[ 
\vert \lambda_{k,2}-\lambda_{k,1}\vert^2\geq C_\varepsilon\lambda_{k,1}^2e^{-2\lambda_{k,1}(T^*+\varepsilon)}.
\]
Finally using this inequality in the estimate of $\sigma_k$  and using 
\[
\frac 1{4T\lambda_{k,i}}\leq \eta(-2T\lambda_{k,i})\leq \frac 1{2T\lambda_{k,i}}
\] 
yield
\[
\sigma_k\geq  \frac{C_\varepsilon}{\lambda_{k,1} + \lambda_{k,2}} e^{-2\lambda_{k,1}(T^*+\varepsilon)}.
\]
From~\eqref{normeminimale},  we deduce
\[
\Vert u_k\Vert^2_{L^2(0,T)}\leq  C_\varepsilon (\lambda_{k,1} + \lambda_{k,2}) e^{-2\lambda_{k,1}(T-(T^*+\varepsilon))} \Vert y_k^0 \Vert^2
\]
which ends the proof.

\end{proof}

%---------------------------------------------------------------------------------------
\subsection{Comments on the behaviour of the spectrum}
\label{Subsec:Oscillateur}
%---------------------------------------------------------------------------------------

We recall some results on the harmonic oscillator proved in~\cite[Sec.~5.1]{DuyckaertsMiller12}. 

Let $H = L^2(\R)$ and $U=L^2(\R)$. Define the operator
\[
A = -\partial_{xx} + x^2, \qquad 
D(A) = \{ y \in H^2(\R) ; x \mapsto x^2 y(x) \in L^2(\R)\}.
\]
The operator is self-adjoint and generates a $C_0$ semigroup. Its eigenvalue are $\{\lambda_k=2k-1 ; k \in \N^*\}$ and its eigenfunctions are the Hermite polynomials (see for instance~\cite{Szego}). Thus they form a Hilbert basis.

Consider the bounded operator $B$ defined by $B = \mathbf{1}_{(-\infty,x_0)}$ with $x_0 \in \R$. 
In this setting it is proved in~\cite[Proposition 5.1]{DuyckaertsMiller12} that system~\eqref{SystControlAbstrait} is not null controllable for any time $T>0$. 

However, still in~\cite[Proposition 5.1]{DuyckaertsMiller12} it is proved that the following inequality holds: there exists $M, m \in \R$ such that
\[
\|v\|^2 \leq M \|(A^* + \lambda) v\|^2 + m \|B^* v \|_{U}^2, \qquad \forall v \in D(A), \: \forall \lambda \in \R.
\]
As for any $T>0$ the functions 
\[
\lambda \in (0,+\infty) \mapsto \frac{e^{2\lambda T}}{\lambda^2} \quad \text{and} \quad
\lambda \in (0,+\infty) \mapsto \frac{e^{2\lambda T}}{\lambda}
\]
are bounded from below,  following definition~\eqref{DefTempsMinimal}, it comes that $T^*=0$. Thus inequality~\eqref{InegResolvante} is valid for any $T$ whereas null controllability holds for no value of $T$. 

It is fundamental to notice that this example satisfies every assumption of Theorem~\ref{ThCondensationGeneral} except that in this situation 
\[
\sum_{k \in \N^*} \frac{1}{\lambda_k} = \sum_{k \in \N^*} \frac{1}{2k-1} = +\infty.
\]
Thus assumption~\eqref{ComportementSpectre} should not be considered as a technical assumption due to the use of moment method but as a concrete limitation of the characterization of null controllability using inequality~\eqref{InegResolvante}.

%%%%%%%%%%%%%%%%%%%%%%%%%%%%%%%%%%%%%%%%%%%%%%%%%%%%%%%%%%%%%%%%%%%%%%%%%%%%%%%%%%%%%%%%%%%%%%%%%%%%%%%%%%%%%%%%%%%%%%%%%%%%%%%%%
%%%%%%%%%%%%%%%%%%%%%%%%%%%%%%%%%%%%%%%%%%%%%%%%%%%%%%%%%%%%%%%%%%%%%%%%%%%%%%%%%%%%%%%%%%%%%%%%%%%%%%%%%%%%%%%%%%%%%%%%%%%%%%%%%
%%%%%%%%%%%%%%%%%%%%%%%%%%%%%%%%%%%%%%%%%%%%%%%%%%%%%%%%%%%%%%%%%%%%%%%%%%%%%%%%%%%%%%%%%%%%%%%%%%%%%%%%%%%%%%%%%%%%%%%%%%%%%%%%%
\appendix
\section{General estimate of biorthogonal families}
\label{AnnexeBiortho}

Let us start recalling that the Blaschke product associated with the sequence $\Lambda = (\lambda_k)_{k \in \mathbb{N}^*}$ is the function $W:\mathbb{C}_{+} \rightarrow \mathbb{C}$ defined by:
	\begin{equation*}
	\left\{ 
	\begin{array}{l}
\displaystyle W(\lambda )=W(\lambda ,\Lambda )= \prod_{k\geq 1} \delta_k  \frac{1-\lambda /\lambda_{k}}{1+\lambda /\overline \lambda_{k}}, \quad \lambda \in  \mathbb{C}_+\, , \\ 
	\noalign{\smallskip} 
\delta_k = \displaystyle \frac{\lambda_k}{\overline \lambda_k} \frac{\left| \lambda_k - 1 \right|}{\left| \lambda_k + 1 \right|}\frac{\overline \lambda_k + 1}{\overline \lambda_k - 1 } \quad (\delta_k = 1 \mbox{ if } \lambda_k = 1).
	\end{array}
	\right.
	\end{equation*}
Under assumption~\eqref{ComportementSpectre}, the previous function $W$ is well-defined, $W \in H^\infty (\mathbb{C}_+)$, the space of bounded and holomorphic functions defined on $\mathbb{C}_+$, is defined almost everywhere on $i \mathbb{R}$ and satisfies $\left\vert W(\lambda )\right\vert <1$, for $\Re {\lambda }>0$, and $\left\vert W(i\tau )\right\vert =1$, for almost every $\tau \in \mathbb{R}$ (see for instance~\cite[Lemma 4.2]{AKBGBDT_Kalman}). A direct computation also gives
	$$
| W' (\lambda_k) | = \frac{1}{2 \mathrm{Re} (\lambda_k)} \mathcal{P}_k^{-1}, \quad \hbox{with } \mathcal{P}_k:= \prod_{ \substack{ l \ge 1  \\ l \not= k }} \left|  \frac{1 + \lambda_k / \overline \lambda_l }{1 - \lambda_k /\lambda_l } \right| . 
	$$
On the other hand (see~\cite[Theorem 3.9]{AKBGBDT_condensation}),
	\begin{equation}\label{f1}
c ( \Lambda ) = \limsup_{k \to + \infty} \frac{- \ln | W ' (\lambda_k)|}{\mathrm{Re} (\lambda_k)} = \limsup_{k \to + \infty} \frac{\ln \mathcal{P}_k }{\mathrm{Re} (\lambda_k)}.
	\end{equation}

Let us now prove Proposition~\ref{Prop:FamilleBiorthogonaleJordan} when $T= \infty$. To this end, let us first consider functions $e_{k,1} = e^{- \lambda_k t}$ and $e_{k,2} = t e^{- \lambda_k t }$ in $L^2(0, \infty ; \mathbb C)$.  Assumption~\eqref{ComportementSpectre} allows us to apply~\cite[Proposition 4.1]{AKBGBDT_Kalman} (with $\eta =2$), and deduce the existence of a biorthogonal family $\{ \widetilde q_{k,1}, \widetilde q_{k,2} \}_{k \in \mathbb{N}^*} \subset A(\Lambda, \infty)$ to $\{ e_{k,1}, e_{k,2} \}_{k \in \mathbb{N}^*}$ which satisfies for any $k \in \N^*$ and any $j \in \{1,2\}$
	\begin{equation}\label{f2}
\| \widetilde q_{k,j} \|_{L^2(0,\infty; \mathbb{C})} \le C  \left[ 1 + \frac{1}{\mathrm{Re} (\lambda_k)} \right] \left(\mathrm{Re} (\lambda_k) \right)^{4} |1 + \lambda_k|^{8} \mathcal{P}_k^{4}, 
	\end{equation}
where $C>0$ is a constant and, for $\tau \in (0, \infty]$, $A ( \Lambda , \tau )$ is the closed subspace of $L^{2} ( 0,\tau; \mathbb{C} )$ given by 
	\begin{equation*}
A ( \Lambda , \tau ) = \overline{\mathrm{span} \, \left\{ e_{k,1}, e_{k,2}: k \in \mathbb{N}^* \right\}}^{L^{2} ( 0,\tau ; \mathbb{C} ) }.
	\end{equation*}
From formulae~\eqref{f1} and~\eqref{f2}, we deduce that, for any $\varepsilon >0$, there exists a constant $C_\varepsilon > 0$ such that
	\begin{equation}\label{f3}
\| \widetilde q_{k,j} \|_{L^2(0,\infty; \mathbb{C})} \le C_\varepsilon e^{(4 c(\Lambda ) + \varepsilon ) \mathrm{Re} (\lambda_k) }, \quad \forall k \in \mathbb{N}^*, \quad \forall j \in \{1,2\} .
	\end{equation}
This proves Proposition~\ref{Prop:FamilleBiorthogonaleJordan} and~\eqref{EstimeeBiorthoJordan} when $T= \infty$.

The general case $T \in (0, \infty)$ can be deduced from the previous result and 
\begin{lemme}\label{isomorphism} 
Let $\Lambda = (\lambda_k)_{k \in \mathbb{N}^*}$ be a sequence satisfying~\eqref{ComportementSpectre}. Then, for any $T \in (0, \infty) $, the restriction operator
	\begin{equation*}
R_{T}: \varphi \in A( \Lambda ,\infty ) \longmapsto R_T \varphi =  \varphi _{\mid ( 0,T) } \in A( \Lambda ,T) 
	\end{equation*}
is a topological isomorphism. In particular, there exists a positive constant $C_{T}$, depending on the sequence $\Lambda$ and $T$, such that
	\begin{equation*}
\Vert \varphi \Vert_{L^{2} ( 0,\infty ; \mathbb{C} ) }\leq C_{T} \Vert \varphi \Vert_{L^{2} ( 0,T ; \mathbb{C} ) },\quad \forall \varphi \in A( \Lambda ,\infty ) . \eqno \Box
	\end{equation*}

\end{lemme}

Before proving Lemma~\ref{isomorphism}, let us complete the proof of Proposition~\ref{Prop:FamilleBiorthogonaleJordan}. Let us set 
	\begin{equation*}
q_{k,j} =  \left(R_T^{-1}\right)^* \widetilde q_{k,j} \in A( \Lambda, T) , \quad \forall k \in \mathbb{N}^*, \: \forall j \in \{1, 2\}.
	\end{equation*}
From Lemma~\ref{isomorphism} and~\eqref{f3}, it is clear that the function $q_{k,j}$ satisfies, for any $\varepsilon >0$, inequality~\eqref{EstimeeBiorthoJordan} for any $ k \in \mathbb{N}^*$ and $j \in \{ 1,2\}$.

On the other hand, we can write 
	\begin{equation*}
\begin{aligned}
\delta_{k,l} \delta_{j,i} &= ( e_{k,j} \, , \, \widetilde q_{l,i} )_{L^2(0, \infty ; \mathbb{C} )} 
\\
&= (R_T^{-1} R_Te_{k,j} \, , \, \widetilde  q_{l,i} )_{L^2(0, \infty ; \mathbb{C} )} \\ 
 &= ( R_T e_{k,j} \, , \, \left(R_T^{-1} \right)^* \widetilde q_{l,i} )_{L^2(0, T ; \mathbb{C} )} 
 \\
 &= (e_{k,j} \, , \, q_{l,i} )_{L^2(0, T ; \mathbb{C} )}, \quad \forall  k,  l \in \mathbb{N}^*, \forall j,i \in \{ 1,2\}
\end{aligned}
	\end{equation*}
i.e. $\left\{ q_{k,1}, q_{k,2}\right\}_{k \in \mathbb{N}^* }\subset A (\Lambda, T) $ is a biorthogonal family to $\left\{e_{k,1}, e_{k,2} \right\}_{k \in \mathbb{N}^*}$ in $L^2(0,T; \mathbb{C})$. This ends the proof of Proposition~\ref{Prop:FamilleBiorthogonaleJordan}.

\begin{proof}[Proof of Lemma~\ref{isomorphism}]
The proof follows the same arguments of the proof of~\cite[Lemma 4.2]{AKBGBDT_condensation} and uses \cite[Lemma 4.6]{AKBGBDT_condensation}.

From the definition of the space $A (\Lambda, T)$, it is clear that we will obtain the proof if we show that there exits a positive constant $C_T >  0$ such that
	$$
\| P \|_{L^2 (0, \infty;  \mathbb{C})} \le C_T \| P \|_{L^2 (0, T;  \mathbb{C})}, \quad \forall P \in \mathfrak{P},
	$$
where
	$$
\mathfrak{P}:= \left\{ P : P(z) = \sum_{j=1}^N \left( a_j e^{-\lambda_j z} + b_j z e^{-\lambda_j z} \right), \hbox{ with } N \ge 1, \ a_j, b_j \in \mathbb{C} \right\}.
	$$
We argue by contradiction. Therefore, assume that there exists a sequence $\{ P_m \}_{m \ge 1} \subset \mathfrak P$ such that
	\begin{equation}\label{f4}
\lim \| P_{m}\|_{L^{2} ( 0,T ; \mathbb{C}) } = 0 \quad \hbox{and} \quad \| P_{m}\|_{L^{2} ( 0,\infty ; \mathbb{C} ) } = 1 \quad 
\forall m \geq 1  . 
	\end{equation}

Observe that condition~\eqref{ComportementSpectre} implies the existence of $ \theta_\delta \in [0, \pi/2) $ such that 
	\begin{equation}\label{sector}
\Lambda =\left\{ \lambda_{k} \right\}_{k\geq 1} \subset S_\delta := \{ z = r e^{i \theta} \in \mathbb{C} : r >0, \  |\theta | \le \theta_\delta \}.
	\end{equation}

Let $\theta _{0}\in \left( \theta _{\delta },\pi /2\right) $ where $\theta_{\delta }\in  ( 0,\pi /2 ) $ is such that~\eqref{sector} holds. Let us also fix $ \varepsilon >0$, $\theta_0  \in (\theta_\delta, \pi/2)$ and $\tau \in (0, \cos \theta_0)$ and consider the set
	\begin{equation*}
S_{\varepsilon ,\theta _{0},\tau }=\left\{ z= x+iy : x\geq \varepsilon ,\ \frac{|y | }{x}\leq \frac{\cos \theta _{0}- \tau }{\sin \theta_{0}} \right\}.
	\end{equation*}
Applying \cite[Lemma~4.6]{AKBGBDT_condensation} in the set $S_{\varepsilon ,\theta _{0},\tau }$ we deduce that there exists a constant $C_{\varepsilon }>0$ such that, for any $P \in \mathfrak{P} $, one has
	\begin{equation*}
| P(z) | \leq C_{\varepsilon }\| P \|_{L^{2}( 0,\infty ; \mathbb{C}) }\left( 1 + |z| \right) e^{- \frac14 \alpha \tau \mathrm{Re} (z) }, \quad \forall z \in S_{\varepsilon ,\theta _{0}, \tau},
	\end{equation*}
where $ \alpha = \min_{k\geq 1}| \lambda_k | >0  $. Using now~\eqref{f4}, we can conclude that the sequence $\left\{ P_{m}\right\}_{m \ge 1} $ is uniformly bounded on the domain $S_{\varepsilon ,\theta_0, \tau }$. Therefore, it is a normal family of holomorphic functions on $S_{\varepsilon ,\theta_0, \tau }$ and there exists a subsequence, still denoted by $\left\{ P_{m}\right\}_{m \ge 1}$, and a holomorphic function $P$ on $S_{\varepsilon ,\theta_0, \tau }$ such that $P_{m}\rightarrow P$ uniformly on the compact sets of $S_{\varepsilon ,\theta_0, \tau }$. Furthermore, from Lebesgue's theorem, $ P_{m}\rightarrow P$ in $L^{2} ( \eta ,\infty ; \mathbb{C}) $ for any $\eta > \varepsilon $. Assumption~\eqref{f4} implies that $P\equiv 0$ on the interval $ ( \eta ,T ) $ for any $\eta : 0<\varepsilon <\eta <T$. Since $P$ is holomorphic on $S_{\varepsilon ,\theta_0, \tau }$, we get $P\equiv 0$ on $ ( \varepsilon ,\infty  ) $. Whence $ \lim \| P_{m}\|_{L^{2}  ( T, \infty; \mathbb{C}  ) }=0$ and since, by our assumption, $\lim \| P_{m}\|_{L^{2} ( 0, T ; \mathbb{C} ) }=0$ it follows that 
	\begin{equation*}
\lim \| P_{m}\|_{L^{2} ( 0,\infty  ; \mathbb{C}) }= 0.
	\end{equation*}
This contradicts~\eqref{f4} and provides the proof of Lemma~\ref{isomorphism}. 

\end{proof}

\section{A sequence with zero condensation index but no gap}
\label{AnnexeNoGap}

\begin{prop} \label{lembohr}
Let 
\[
\Lambda=\{k^2, k^2+\alpha_k,\,k\geq 1\}=\{\lambda_{2k-1}:=k^2,\,\lambda_{2k}=k^2+\alpha_k,\,k\geq 1\},
\]
with $\alpha_k\underset{k\to +\infty}{\longrightarrow} 0$ and $\alpha_k\in (0,1)$. 
Then, 
\begin{equation}
\operatorname{c}(\Lambda)= \underset{k \rightarrow +\infty} { \limsup} \,  \frac{-\ln { \vert E'(\lambda_k)  \vert } }{\lambda_k}= \underset{k\to +\infty} \limsup   \frac{-\ln (\alpha_k)}{k^2},
\label{bohr}
\end{equation}
where $E$ is defined by~\eqref{DefE}.
% \begin{equation}  C(\lambda)=\underset {k\geq1} {\prod}  (1-\frac{\lambda^2}{\lambda_k^2}).\label{interpolation} \end{equation}
\end{prop}

This proposition is proved below by direct computations inspired by discussions with L.~\textsc{Ouaili}\footnote{Aix Marseille Universit\'e, CNRS, Centrale Marseille, I2M, UMR 7373}.

\begin{proof}
Direct computations give
\begin{align*}
\frac{\ln \vert  E' ( \lambda_{m})\vert} {\lambda_{m}} 
&= \frac{1}{\lambda_m} \ln \left|\left(\frac{-2}{\lambda_m} \right) \prod_{j \neq m} \left(1- \frac{\lambda_m^2}{\lambda_j^2} \right) \right|
\\
&=  \frac{\ln \vert \lambda_{m}-\lambda_{m+1} \vert}{\lambda_{m}} +  \frac{\ln \vert \lambda_{m}-\lambda_{m-1} \vert}{\lambda_{m}} 
\\
&+ \frac{1}{\lambda_{m}} \ln  \left(\frac{2\left(\lambda_{m}+\lambda_{m+1} \right)\left(\lambda_{m-1}+\lambda_m\right)) }{\lambda_m\lambda^2_{m-1}\lambda^2_{m+1}}\right)+ \,F_m+ \,G_m 
\end{align*}
where 
\[
F_m:=\underset { j < m-1}\sum \frac {\ln\vert 1-\frac{\lambda_{m}^2}{\lambda_{j}^2} \vert}{\lambda_{m}},
\qquad 
G_m:= \underset { j>m+1}\sum \frac {\ln\vert 1-\frac{\lambda_{m}^2}{\lambda_{j}^2} \vert}{\lambda_{m}}.
\]

First, taking into account the expression of $\lambda_m$ for $m $ even or odd, it is clear that 
	$$
\underset{m\to +\infty} \limsup \left(   \frac{\ln \vert \lambda_{m}-\lambda_{m+1} \vert}{\lambda_{m}}+  \frac{\ln \vert \lambda_{m}-\lambda_{m-1} \vert}{\lambda_{m}} \right) 
=  \limsup_{k \to +\infty} \frac{\ln (\alpha_k)}{ k^2}.
	$$

Secondly, from the definition of $\lambda_{m}$, one gets that 
	$$
\underset{m\to +\infty}\limsup \frac{1}{\lambda_{m}} \ln  \left(\frac{2\left(\lambda_{m}+\lambda_{m+1} \right)\left(\lambda_{m-1}+\lambda_m\right)) }{\lambda_m \lambda^2_{m-1}\lambda^2_{m+1}}\right)=0.
	$$

In order to end the proof of the result, let us see that $F_m$ and $G_m$ go to $0$ when $m$ goes to $+\infty$. 

\begin{itemize}
\item[$\bullet$] Study of $F_m$. Notice that
\[
\vert F_m\vert \leq \sum\limits_{2\lambda_j^2 {\le} \lambda_m^2} \frac 1{\lambda_m} \ln\left(\frac {\lambda_m^2}{\lambda_j^2}-1\right) + 
\sum_{2\lambda_j^2>\lambda_m^2,\,j<m-1} \frac 1{\lambda_m}\ln\left(\frac {\lambda_j^2}{\lambda_m^2-\lambda_j^2}\right)
\]
The first term is estimated by
\begin{align*}
\sum\limits_{2\lambda_j^2 {\le} \lambda_m^2} \frac 1{\lambda_m}\ln\left(\frac {\lambda_m^2}{\lambda_j^2}-1\right)
&\leq \sum_{j<m-1} \frac 1{\lambda_m}\ln\left(\frac {\lambda_m^2}{\lambda_1^2}-1\right)=\frac {m-1}{\lambda_m}\ln\left( {\lambda_m^2}-1\right)
\\
&\leq 2(m-1)\frac{\ln(\lambda_m )}{\lambda_m}\underset{m\to +\infty}{\longrightarrow} 0,
\end{align*}
and the second term is estimated by
\begin{align*}
\sum_{2\lambda_j^2>\lambda_m^2,\,j<m-1} \frac 1{\lambda_m}\ln\left(\frac {\lambda_j^2}{\lambda_m^2-\lambda_j^2}\right)
&\leq \underset{j<m-1}\sum \frac 1{\lambda_m}\ln\left(\frac {\lambda_j^2}{\lambda_m^2-\lambda_{m-2}^2}\right)
\\
&\leq \frac{m-1}{\lambda_m}\ln \left(\frac{\lambda_m^2}{\lambda_m^2-\lambda_{m-2}^2}\right)\underset{m\to +\infty}{\longrightarrow} 0.
\end{align*}
Thus, 
\[
\underset{m\to +\infty}\lim F_m=0.
\]

\item[$\bullet$] Study of $G_m$. Notice that
\[
\vert G_m\vert \leq \underset{j>m+1}\sum \frac 1{\lambda_m}\ln\left(\frac {\lambda_j^2}{\lambda_j^2-\lambda_m^2}\right)=\underset{j>m+1}\sum \frac 1{\lambda_m}\ln\left(1+\frac{\lambda_m^2}{\lambda_j^2-\lambda_m^2}\right).
\]
Using the inequality $\ln(1+x) \leq x$, when $x >0$, this reduces to
\begin{align*}
\vert G_m\vert
&\leq \underset{j\geq m+2}\sum \frac{\lambda_m}{\left(\lambda_j-\lambda_m\right)\left(\lambda_j+\lambda_m\right)}
\leq \underset{j\geq m+2}\sum \frac 1{\lambda_j-\lambda_m}
\\
&=\sum_{j= {2p - 1}\geq m+2} \frac 1{ p^2 - \lambda_m } + \sum_{j={2p}\geq m+2} \frac 1{ p^2+\alpha_p-\lambda_m }\\
&  \le  \sum_{p \ge (m+3)/2 } \frac 1{ p^2 - \lambda_m }+ \sum_{p \ge (m+2)/2 } \frac 1{ p^2 + \alpha_p - \lambda_m }\\
& \le 2  \sum_{p \ge (m+2)/2 } \frac 1{ p^2 - \lambda_m } := H_m, \quad \forall k \ge 1.
\end{align*}

It is not difficult to check that
	$$
H_{2k-1}= 2 \sum_{p \ge k+1} \frac 1{p^2 - k^2} ; \quad H_{2k}= 2 \sum_{p \ge k+1} \frac 1{p^2 - k^2- \alpha_k }, \quad \forall k \ge 1.
	$$

Taking into account that $ H_{2k-1} < H_{2k}$ for any $k \ge 1$, let us estimate $H_{2k}$. To this end, let us consider the function $f_k :  [k+1,+\infty) \to \R$ given by
	$$
f_k (t) = \frac 1{t^2 - k^2 - \alpha_k}, \quad \forall t \in [k+1,+\infty).
	$$
Observe that $f_k$ is a nonincreasing function on $ [k+1,+\infty) $. Therefore
	$$
	\begin{array}{l}
\displaystyle 0 \le H_{2k}= 2 \sum_{p \ge k+1} f_k (p) \le 2 \left( f_k(k+1) + \int_{k+1}^{+ \infty} f_k (t) \, \mathrm{d} t \right) \\
	\noalign{\smallskip}
\displaystyle \phantom{\displaystyle 0 \le H_{2k}} = \frac 2{2k+1 - \alpha_k} + \frac 1{\sqrt{k^2 + \alpha_k}} \ln \left( \frac{(k+1+ \sqrt{k^2 + \alpha_k} )^2}{2k+1 - \alpha_k} \right). 
	\end{array}
	$$

From the previous computations we can conclude 
	\[
\underset{m\to +\infty}\lim G_m=0.
	\]

\end{itemize}

In conclusion, we have proved that 
	\[
\operatorname{c}(\Lambda)= \underset{k\to +\infty} \limsup   \frac{-\ln (\alpha_k)}{k^2}.
	\]
This finalizes the proof of Proposition~\ref{lembohr}. 

\end{proof}

\section{Condensation and Bohr index}
\label{AnnexeBohr}

Consider the sequence defined by $\Lambda = \{ \lambda_k, d \lambda_k ; k \in \N^*\}$ with $\lambda_k = k^2 \pi^2$ and $\sqrt{d} \not\in \mathbb{Q}$. 

For the sake of completeness, we reproduce here the computations of~\cite{Samb17} communicated by E.H.~\textsc{Samb} which proves that $\operatorname{c}(\Lambda) = \operatorname{Bohr}(\Lambda)$. 

From the computations in the proof of~\cite[Proposition 6.20]{AKBGBDT_condensation}, it comes that 
\[
\operatorname{c}(\Lambda)= \max\{ \ell_1, \ell_2 \},
\]
with 
\[
\ell_1 = \limsup\limits_{k \to +\infty} \frac{-\ln |E'(\lambda_k)|}{\lambda_k}, \qquad
\ell_2 = \limsup\limits_{k \to +\infty} \frac{-\ln |E'(d\lambda_k)|}{d\lambda_k}.
\]
It is also proved that
\[
\ell_1 = \limsup\limits_{k \to +\infty} \frac{-\ln \left| \sin \left( \frac{k \pi}{\sqrt{d}} \right) \right|}{\lambda_k}, \qquad
\ell_2 = \limsup\limits_{k \to +\infty} \frac{-\ln \left| \sin \left( k \pi \sqrt{d} \right) \right|}{d\lambda_k}.
\]

The fact that the condensation index is indeed a Bohr index follows from the application of the next lemma to the particular sequence $(k_n)_{n \in \N^*}$ such that
\[
\lim\limits_{n \to +\infty} \frac{-\ln |E'(\lambda_{k_n})|}{\lambda_{k_n}} = \ell_1.
\]
\begin{lemme}[\cite{Samb17}] \label{Lem:Cond=Bohr}
For any sequence of integers $(k_n)_{n \in \N^*}$ going to $+\infty$ as $n$ goes to $+\infty$ there exists a sequence $(j_n)_{n \in \N^*}$ such that
\[
\lim\limits_{n \to +\infty} \left| \frac{-\ln |\lambda_{k_n}-d \lambda_{j_n}|}{\lambda_{k_n}} - \frac{-\ln |E'(\lambda_{k_n})|}{\lambda_{k_n}} \right| = 0.
\]
\end{lemme}

\begin{proof}
For any $n \in \N^*$, let $j_n$ the nearest integer to $\frac{k_n}{\sqrt{d}}$. Since 
\[
\left| \frac{k_n \pi}{\sqrt{d}} - j_n \pi \right| \leq \frac{\pi}{2},
\]
we obtain
\[
2 \left| \frac{k_n}{\sqrt{d}} -j_n \right| \leq \sin \left| \frac{k_n \pi}{\sqrt{d}} - j_n \pi \right| =  \left| \sin \left( \frac{k_n \pi}{\sqrt{d}} \right) \right| \leq \pi \left| \frac{k_n}{\sqrt{d}} - j_n \right|.
\]
Using~\cite[Proposition 6.20]{AKBGBDT_condensation} leads to 
\[
|E'(\lambda_k)| = \frac{d \sinh(k \pi) \sinh \left( \frac{k \pi}{\sqrt{d}} \right)}{2 k^5 \pi^3} \left| \sin\left( \frac{k \pi}{\sqrt{d}} \right) \right|, \quad \forall k \in \N^*.
\]
Finally, this yields
\begin{align*}
&\left| \frac{-\ln |\lambda_{k_n}-d \lambda_{j_n}|}{\lambda_{k_n}} - \frac{-\ln |E'(\lambda_{k_n})|}{\lambda_{k_n}} \right|
= \frac{1}{\lambda_{k_n}} \left| \ln \left| \frac{E'(\lambda_{k_n})}{\lambda_{k_n}-d \lambda_{j_n}} \right|  \right|
\\
&= \frac{1}{\lambda_{k_n}} \left| \ln \left( \frac{d \sinh(k \pi) \sinh \left( \frac{k \pi}{\sqrt{d}} \right)}{2 k^5 \pi^5 (k_n + \sqrt{d} j_n)} \left| \frac{\sin\left( \frac{k \pi}{\sqrt{d}}\right)}{k_n - \sqrt{d} j_n}  \right| \right) \right|
\\
&\leq \frac{1}{\lambda_{k_n}} \left| \ln \left( \frac{d \sinh(k \pi) \sinh \left( \frac{k \pi}{\sqrt{d}} \right)}{2 k^5 \pi^5 (k_n + \sqrt{d} j_n)}  \right) \right| + \frac{1}{\lambda_{k_n}} \left| \ln \frac{\pi}{\sqrt{d}} \right|
\\
&\leq \frac{1}{\lambda_{k_n}} \left| \ln \frac{d e^{k_n \pi} e^{\frac{k_n \pi}{\sqrt{d}}}}{8 k_n^5 \pi^5 (k_n + \sqrt{d} j_n)} \right|+ \frac{1}{\lambda_{k_n}} \left| \ln \frac{\pi}{\sqrt{d}} \right|
\end{align*}
The right-hand side goes to $0$ as $n$ goes to $+\infty$ which ends the proof.

\end{proof}

%\bibliographystyle{plain} 
%\bibliography{/home/morgan/Documents/Recherche/biblio}

\end{document}